\documentclass{article}

\usepackage{fancyhdr, fullpage}
\usepackage[american]{babel}
\usepackage{url}            % simple URL typesetting
\usepackage{booktabs}       % professional-quality tables
\usepackage{amsfonts}       % blackboard math symbols
\usepackage{nicefrac}       % compact symbols for 1/2, etc.
\usepackage{microtype}      % microtypography
\usepackage{natbib}

\usepackage[breaklinks,colorlinks,citecolor=blue]{hyperref}
\usepackage{amssymb, amsmath, amsthm}
\usepackage{algorithm}
\usepackage{algpseudocode}
\usepackage{mathtools}
\usepackage{enumitem}
\usepackage{colortbl}
\usepackage{xcolor}
\usepackage{ragged2e}
\usepackage{multirow}
\usepackage{array, hhline}
\usepackage{subcaption, graphicx}
\usepackage{siunitx}

\newcommand{\CC}{C\nolinebreak\hspace{-.05em}\raisebox{.4ex}{\tiny\bf +}\nolinebreak\hspace{-.10em}\raisebox{.4ex}{\tiny\bf +}}
\newcommand{\R}{\mathbb{R}}

\newcommand{\pf}[1]{\nabla f (#1)}

\newcommand{\pfik}[1]{\nabla f_{i_k}(#1)}

\newcommand{\pfi}[1]{\nabla f_{i}(#1)}

\newcommand{\norm}[1]{\left\lVert#1\right\rVert}

\newcommand{\innr}[1]{\left\langle#1\right\rangle}

\newcommand{\Eik}[1]{\mathbb{E}_{i_k} \left[#1\right]}
\newcommand{\Ei}[1]{\mathbb{E}_{i} \left[#1\right]}

\newcommand{\Er}[2]{\mathbb{E}_{#1} \big[#2\big]}
\newcommand{\E}[1]{\mathbb{E} \left[#1\right]}
\newcommand{\mleq}[1]{\overset{\mathclap{(#1)}}{\leq}}
\newcommand{\meq}[1]{\overset{\mathclap{(#1)}}{=}}

\newcommand{\mra}[1]{\overset{\mathclap{(#1)}}{\Rightarrow}}
\newcommand{\mar}[1]{\left(#1\right)}
\newcommand{\xs}{x^{\star}}

\newcommand{\supp}[1]{[#1]_{T_{i_k}}}
\newcommand{\suppi}[1]{[#1]_{T_i}}

\graphicspath{{./Figs/}}

\newcommand{\Title}{Accelerating Perturbed Stochastic Iterates in Asynchronous Lock-Free Optimization}

\newtheorem{lemmas}{Lemma}
\newtheorem{theorems}{Theorem}

\newtheorem{pro-remark}{Remark}[propositions]
\newtheorem{corollarys}{Corollary}[theorems]
\newtheorem{assumptions}{Assumption}

\begin{document}

\title{\Title}

\author{Kaiwen Zhou\thanks{Department of Computer Science and Engineering, The Chinese University of Hong Kong, Sha Tin, N.T., Hong Kong SAR; e-mail: \href{mailto:kwzhou@cse.cuhk.edu.hk}{\tt kwzhou@cse.cuhk.edu.hk}.} \and Anthony Man-Cho So\thanks{Department of Systems Engineering and Engineering Management, The Chinese University of Hong Kong, Sha Tin, N.T., Hong Kong SAR; e-mail: \href{mailto:manchoso@se.cuhk.edu.hk}{\tt manchoso@se.cuhk.edu.hk}.} \and James Cheng\thanks{Department of Computer Science and Engineering, The Chinese University of Hong Kong, Sha Tin, N.T., Hong Kong SAR; e-mail: \href{mailto:jcheng@cse.cuhk.edu.hk}{\tt jcheng@cse.cuhk.edu.hk}.}}

\pagestyle{plain}
\date{}

\maketitle

\begin{abstract}
	We show that stochastic acceleration can be achieved under the perturbed iterate framework \citep{ManiaPPRRJ17} in asynchronous lock-free optimization, which leads to the optimal incremental gradient complexity for finite-sum objectives. We prove that our new accelerated method requires the same linear speed-up condition as the existing non-accelerated methods. Our core algorithmic discovery is a new accelerated SVRG variant with sparse updates. Empirical results are presented to verify our theoretical findings. 
\end{abstract}

\section{Introduction}
We consider the following unconstrained optimization problem, which appears frequently in machine learning and statistics, such as empirical risk minimization:
\begin{equation}\label{prob_def}
	\min_{x\in \R^d}{f(x) \triangleq \frac{1}{n} \sum_{i=1}^n{f_i(x)}},
\end{equation}
where each $f_i$ is $L$-smooth and convex, $f$ is $\mu$-strongly convex.\footnote{In fact, we will only use a weaker quadratic growth assumption. The formal definitions are in Section \ref{sec:preliminaries}. } We further assume that the computation of each $f_i$ is \textit{sparse}, i.e., the computation is supported on a set of coordinates $T_i \subseteq \{1,\ldots, d\}$. For example, generalized linear models (GLMs) satisfy this assumption. GLMs has the form: $\forall i, f_i(x) = \ell_i(\innr{a_i, x})$, where $\ell_i$ is some loss function and $a_i\in \R^d$ is a data sample. In this case, the support of each $f_i$ is the set of non-zero coordinates of $a_i$. Many large real-world datasets have very high sparsity (see Table \ref{data-table} for several examples).

In the past decade, exciting progress has been made on solving problem \eqref{prob_def} with the introduction of stochastic variance reduced methods, such as SAG \citep{SAG,SAG2}, SVRG \citep{SVRG,Prox_SVRG}, SAGA \citep{SAGA}, S2GD \citep{S2GD}, SARAH \citep{SARAH}, etc. These methods achieve fast linear convergence rate as gradient descent (GD), while having low iteration cost as stochastic gradient descent (SGD). There are also several recent works on further improving the practical performance of these methods based on past gradient information \citep{HofmannLLM15,TanMDQ16,ShangLZCNY18,shi2021ai,dubois2021svrg}. Inspired by Nesterov's seminal works \citep{AGD1,AGD2,AGD3}, stochastic accelerated methods have been proposed these years \citep{allen2017katyusha,MiG,SSNM,pmlr-v124-zhou20a,VARAG,VRADA,Anita}. These methods achieve the optimal convergence rates that match the lower complexity bounds established in \citep{woodworthS16}. For sparse problems, lagged update \citep{SAG2} is a common technique for the above methods to efficiently leverage the sparsity.

Inspired by the emerging parallel computing architectures such as multi-core computer and distributed system, many parallel variants of the aforementioned methods have been proposed, and there is a vast literature on those attempts. 
%The textbook \citep{bertsekas1989parallel} provides the foundational study for parallel algorithms. 
Among them, asynchronous lock-free algorithms are of special interest.
\citet{Hogwild} proposed an asynchronous lock-free variant of SGD called Hogwild!\ and first proved that it can achieve a \textit{linear speed-up}, i.e., running Hogwild!\ with $\tau$ parallel processes only requires $O(1 / \tau)$-times the running time of the serial variant. The condition on the maximum overlaps among the parallel processes is called the \textit{linear speed-up condition}. Following \citet{Hogwild}, \citet{SaZORR15,LianHLL15} analyzed asynchronous SGD for non-convex problems, and \citet{DuchiJM13} analyzed for stochastic optimization; \citet{ManiaPPRRJ17} refined the analysis framework (called the perturbed iterate framework) and proposed KroMagnon (asynchronous SVRG); \citet{ASAGA} simplified the perturbed iterate analysis and proposed ASAGA (asynchronous SAGA), and in the journal version of this work,  \citet{ASAGA_JMLR} further improved the analysis of KroMagnon and Hogwild!; \citet{ProxASAGA} derived proximal variant of ASAGA; \citet{NguyenNDRST18} refined the analysis of Hogwild!; \citet{JoulaniGS19} proposed asynchronous methods in online and stochastic settings; \citet{gu2020unified} proposed several asynchronous variance reduced algorithms for non-smooth or non-convex problems; \citet{StichMJ21} studied a broad variety of SGD variants and improved the speed-up condition of asynchronous SGD.

All the above mentioned asynchronous methods are based on non-accelerated schemes such as SGD, SVRG and SAGA. Witnessing the success of recently proposed stochastic accelerated methods, it is natural to ask: \textit{Is it possible to achieve stochastic acceleration in asynchronous lock-free optimization? Will it lead to a worse linear speed-up condition?} The second question comes from a common perspective that accelerated methods are less tolerant to the gradient noise, e.g., \citep{devolder2014first,RMM}.

In this work, we give positive answers to these questions, i.e., we propose the first asynchronous stochastic accelerated method for solving problem~\eqref{prob_def}, and we prove that it requires the identical linear speed-up condition as the non-accelerated counterparts. 

The works that are closest related to us are \citep{fang2018accelerating,MiG,A2BCD}.  \citet{fang2018accelerating} proposed several asynchronous accelerated schemes. However, their analysis requires consistent read and does not consider sparsity. The dependence on the number of overlaps $\tau$ in their complexities is $O(\tau)$ compared with $O(\sqrt{\tau})$ in our result. \citet{MiG} naively extended their proposed accelerated method into asynchronous lock-free setting. However, even in the serial case, their result imposes strong assumptions on the sparsity. This issue is discussed in detail in Section \ref{sec:svc}. In the asynchronous case, no theoretical acceleration is achieved in \citep{MiG}. \citet{A2BCD} proposed the first asynchronous accelerated block-coordinate descent method (BCD). However, as pointed out in Appendix F in \citep{ProxASAGA}, BCD methods perform poorly in sparse problems since they  focus only on a single coordinate of the gradient information.

\section{Preliminaries}
\label{sec:preliminaries}
\paragraph{Notations} We use $\innr{\cdot,\cdot}$ and $\norm{\cdot}$ to denote the standard inner product and the Euclidean norm, respectively. We let $[n]$ denote the set $\{1,2,\ldots,n\}$, $\mathbb{E}$ denote the total expectation and $\mathbb{E}_{i}$ denote the expectation with respect to a random sample $i$. We use $[x]_v$ to denote the coordinate $v$ of the vector $x\in \R^d$, and for a subset $T \subseteq [d]$, $[x]_T$ denotes the collection of $[x]_v$ for $v \in T$. We let $\lceil \cdot \rceil$ denote the ceiling function, $I_d$ denote the identity matrix and $(\alpha)_+ \triangleq \max{\{\alpha, 0\}}$. 

We say that a function $f: \R^d \rightarrow \R$ is \textit{$L$-smooth} if it has $L$-Lipschitz continuous gradients, i.e., $\forall x, y\in \R^d$,
\[
\norm{\pf{x} - \pf{y}} \leq L\norm{x - y}.
\]
An important consequence of $f$ being $L$-smooth and convex is that $\forall x, y\in \R^d$,
\[
f(x) - f(y) - \innr{\pf{y}, x - y}  \geq{} \frac{1}{2L} \norm{\pf{x} - \pf{y}}^2.
\]
We call it the \textit{interpolation condition} following \citep{taylor2017smooth}. A continuously differentiable $f$ is called \textit{$\mu$-strongly convex} if $\forall x, y\in \R^d$,
\[
f(x) - f(y) - \innr{\pf{y}, x - y} \geq \frac{\mu}{2} \norm{x - y}^2.
\]
In particular, the strong convexity at any $x\in \R^d$ and $\xs$ is called the \textit{quadratic growth}, i.e., 
\[
f(x) - f(\xs) \geq \frac{\mu}{2} \norm{x - \xs}^2.
\]
All our results hold under this weaker condition. We define $\kappa \triangleq \frac{L}{\mu}$, which is always called the condition ratio. Other equivalent definitions of these assumptions can be found in \citep{AGD3, karimi2016linear}. 

Oracle complexity refers to the number of incremental gradient computations needed to find an $\epsilon$-accurate solution (in expectation), i.e., $\E{f(x)} - f(\xs) \leq \epsilon$.

\section{Serial Sparse Accelerated SVRG}

We first introduce a new accelerated SVRG variant with sparse updates in the serial case (Algorithm~\ref{alg:SS-Acc-SVRG}), which serves as the base algorithm for our asynchronous scheme. Our technique is built upon the following sparse approximated SVRG gradient estimator proposed in \citep{ManiaPPRRJ17}: For uniformly random $i\in [n]$ and $y, \tilde{x} \in \R^d$, define
\begin{equation}\label{S-SVRG-estimator}
	\mathcal{G}_y \triangleq \pfi{y} - \pfi{\tilde{x}} + D_i\pf{\tilde{x}},
\end{equation}
where $D_i \triangleq P_i D$ with $P_i\in \R^{d\times d}$ being the diagonal projection matrix of $T_i$ (the support of $f_i$) and $D = (\frac{1}{n}\sum_{i=1}^n{P_i})^{-1}$, which can be computed and stored in the first pass. Note that we assume the coordinates with zero cardinality have been removed in the objective function,\footnote{Clearly, the objective value is not supported on those coordinates.} and thus $nI_d \succeq D \succeq I_d$. A diagonal element of $D$ corresponds to the inverse probability of the coordinate belonging to a uniformly sampled support $T_i$. It is easy to verify that this construction ensures the unbiasedness $\Er{i}{D_i\pf{\tilde{x}}} = \pf{\tilde{x}}$. 

This section is organized as follows: In Section \ref{sec:svc}, we summarize the key technical novelty and present the convergence result; in Section \ref{sec:remarks}, we remark on the design of Algorithm \ref{alg:SS-Acc-SVRG}; in Section \ref{sec:general_convex}, we briefly discuss its general convex variant.

\subsection{Sparse Variance Correction}
\label{sec:svc}

\begin{algorithm*}[t]
	\caption{SS-Acc-SVRG: Serial Sparse Accelerated SVRG}
	\label{alg:SS-Acc-SVRG}
	\renewcommand{\algorithmicrequire}{\textbf{Input:}}
	\renewcommand{\algorithmicensure}{\textbf{Initialize:}}
	\begin{algorithmic}[1]
		\Require initial guess $x_0 \in \R^d$, constant $\omega > 1$ which controls the restart frequency.
		\Ensure set the scalars $m, \vartheta, \varphi, \eta,  S, R$ according to Theorem \ref{thm:SS-Acc-SVRG}, initialize the diagonal matrix $D$.
		\For{$r=0, \ldots, R-1$}\Comment{performing restarts}
		\State $\tilde{x}_0 = z^0_0 = x_r$.
		\For{$s=0, \ldots, S-1$}
		\State Compute and store $\pf{\tilde{x}_s}$. 
		\For{$k = 0, \ldots, m-1$}
		\State Sample $i_k$ uniformly in $[n]$ and let $T_{i_k}$ be the support of $f_{i_k}$.
		\State\label{alg-step:coupling} $\supp{y_{k}} = \vartheta\cdot \supp{z^s_k} + \left(1 - \vartheta\right) \cdot\supp{\tilde{x}_s} - \varphi\cdot\supp{D\pf{\tilde{x}_s}}.$ \Comment{sparse coupling}
		\State\label{alg-step:mirror_c} $\supp{z^s_{k+1}} = \supp{z^s_k} - \eta \cdot\big(\pfik{\supp{y_k}} - \pfik{\supp{\tilde{x}_s}} + D_{i_k}\pf{\tilde{x}_s}\big)$.\Comment{sparse update} 
		\EndFor
		\State $\tilde{x}_{s+1} = y_t$ for uniformly random $t\in \{0, 1, \ldots, m-1\}$.
		\State $z^{s+1}_0 = z^s_m$.
		\EndFor
		\State $x_{r+1} = \frac{1}{S}\sum_{s = 0}^{S-1}{\tilde{x}_{s+1}}$.
		\EndFor
		\renewcommand{\algorithmicensure}{\textbf{Output:}}
		\Ensure $x_R$.
	\end{algorithmic}
\end{algorithm*}

Intuitively, the sparse SVRG estimator \eqref{S-SVRG-estimator} will leads to a larger variance compared with the dense one (when $D=I_d$). In the previous attempt, \citet{MiG} naively extends an accelerated SVRG variant into the sparse setting, which results in a fairly strong restriction on the sparsity as admitted by the authors.\footnote{It can be shown that when $\kappa$ is large, almost no sparsity is allowed in their result.} By contrast, sparse variants of SVRG and SAGA in \citep{ManiaPPRRJ17,ASAGA} require no assumption on the sparsity, and they achieve the same oracle complexities as their original dense versions.

Analytically speaking, the only difference of adopting the sparse estimator \eqref{S-SVRG-estimator} is on the variance bound. The analysis of non-accelerated dense SVRG typically uses the following variance bound ($D=I_d$) \citep{SVRG, Prox_SVRG}:
\[
\E{\norm{\mathcal{G}_y - \pf{y}}^2} \leq 4L\big(f(y) - f(\xs) + f(\tilde{x}) - f(\xs)\big).
\]
It is shown in \citep{ManiaPPRRJ17} that the sparse estimator \eqref{S-SVRG-estimator} admits the same variance bound for any $D$. Thus, the analysis in the dense case can be directly applied to the sparse variant, which leads to a convergence guarantee that is independent of the sparsity.

However, things are not as smooth in the accelerated case. To (directly) accelerate SVRG, we typically uses a much tighter variance bound in the form ($D=I_d$) \citep{allen2017katyusha}:
\[
\E{\norm{\mathcal{G}_y - \pf{y}}^2} \!\leq 2L\big(f(\tilde{x}) - f(y) - \innr{\pf{y}, \tilde{x} - y}\big).
\]
Unfortunately, in the sparse case ($D\neq I_d$), we do not have an identical variance bound as before. The variance of the sparse estimator \eqref{S-SVRG-estimator} can be bounded as follows. The proof is given in Appendix \ref{app:lemma_vb}. 
\begin{lemmas}[Variance bound for accelerated SVRG]
	\label{lem:S-SVRG-VB}
	The variance of $\mathcal{G}_y$ (defined at \eqref{S-SVRG-estimator}) can be bounded as
	\begin{equation}\label{S-SVRG-VB}
		\begin{aligned}
			\Ei{\norm{\mathcal{G}_y - \pf{y}}^2}
			\leq{}& 2L\big(f(\tilde{x}) - f(y) - \innr{\pf{y}, \tilde{x} - y}\big)  \!-\! \norm{\pf{y}}^2 \\ 
			&+ \underbrace{2\innr{\pf{y}, D \pf{\tilde{x}}}}_{\mathcal{R}_1} -\innr{\pf{\tilde{x}},D\pf{\tilde{x}}}.
		\end{aligned}
	\end{equation}
\end{lemmas}
In general, except for the dense case, where we can drop the last three terms above by completing the square, this upper bound will always be correlated with the sparsity (i.e., $D$). This correlation causes the strong sparsity assumption in \citep{MiG}. We may consider a more specific case where $D = nI_d$. In this case, the last three terms above can be written as $(n-1)\norm{\pf{y}}^2 - n\norm{\pf{y} - \pf{\tilde{x}}}^2$, which is not always non-positive.

Inspecting \eqref{S-SVRG-VB}, we see that it is basically the term $\mathcal{R}_1$, which could be positive, that causes the issue. We thus propose a novel \textit{sparse variance correction} for accelerated SVRG, which is designed to perfectly cancel $\mathcal{R}_1$. The correction is added to the coupling step (Step \ref{alg-step:coupling}):
\[
y_{k} = \vartheta \cdot z_k + \underbrace{\left(1 - \vartheta\right)\cdot \tilde{x}_s}_{\text{Negative Momentum}} - \underbrace{\varphi \cdot D\pf{\tilde{x}_s}}_{\text{Sparse Variance Correction}}.
\]
This correction neutralizes all the negative effect of the sparsity, in a similar way as how the negative momentum cancels the term $\innr{\pf{y}, \tilde{x} - y}$ in the analysis \citep{allen2017katyusha}. We can understand this correction as a variance reducer that controls the additional sparse variance. Then we have the following sparsity-independent convergence result for Algorithm~\ref{alg:SS-Acc-SVRG}, and its proof is given in Appendix \ref{app:thm_serial}. 

\begin{theorems}
	\label{thm:SS-Acc-SVRG}
	For any constant $\omega > 1$, let $m=\Theta(n), \vartheta=\frac{\sqrt{m}}{\sqrt{\kappa} + \sqrt{m}}, \varphi = \frac{1 - \vartheta}{L}$, $\eta = \frac{1 - \vartheta}{L\vartheta}$ and $S= \left\lceil 2\omega\sqrt{\frac{\kappa}{m}}\right\rceil$. For any accuracy \mbox{$\epsilon > 0$,} we restart $R= O\big(\log{\frac{f(x_0) - f(\xs)}{\epsilon}}\big)$ rounds. Then, Algorithm \ref{alg:SS-Acc-SVRG} outputs $x_R$ satisfying $\E{f(x_R)} - f(\xs) \leq \epsilon$ in oracle complexity
	\[
	O\left(\max{\left\{n, \sqrt{\kappa n}\right\}}\log{\frac{f(x_0) - f(\xs)}{\epsilon}}\right).
	\]
\end{theorems}

This complexity matches the lower bound established in \citep{woodworthS16} (up to a log factor), and is substantially faster than the $O\big((n+\kappa)\log{\frac{1}{\epsilon}}\big)$ complexity of sparse approximated SVRG and SAGA derived in \citep{ASAGA, ASAGA_JMLR} in the ill-conditioned regime ($\kappa \gg n$).

\subsection{Other Remarks about Algorithm \ref{alg:SS-Acc-SVRG}}
\label{sec:remarks}

\begin{itemize}[leftmargin=14pt]
	\item We adopt a restart framework to handle the strong convexity, which is also used in \citep{MiG} and is suggested in \citep{allen2017katyusha} (footnote 9). This framework allows us to relax the strong convexity assumption to quadratic growth. Other techniques for handling the strong convexity such as (i) assuming a strongly convex regularizer and using proximal update \citep{nesterov2005smooth,allen2017katyusha}, and (ii) the direct constructions in \citep{VARAG,zhou2020boosting, Anita} fail to keep the inner iteration sparse, and we are currently not sure how to modify them. 
	
	\item The restarting point $x_{r+1}$ is chosen as the averaged point instead of a uniformly random one because large deviations are observed in the loss curve if a random restarting point is used. 
	
	\item We can also use the correction to fix the sparsity issue in \citep{MiG}, which leads to a slightly different method and somewhat longer proof.  
	
	\item Algorithm \ref{alg:SS-Acc-SVRG} degenerates to a strongly convex variant of Acc-SVRG-G proposed in \citep{zhou2021practical} in the dense case ($D=I_d$), which summarizes our original inspiration. Note that Acc-SVRG-G is designed for minimizing the gradient norm and function value at the same time, which has nothing to do with handling sparsity. Its correction is not necessary if we only want to minimize the function value at the optimal rate, unlike our correction. 
\end{itemize}
\subsection{General Convex Variant}
\label{sec:general_convex}
A general convex ($\mu=0$) variant of Algorithm~\ref{alg:SS-Acc-SVRG} can be derived by removing the restarts and adopting a variable parameter choice similar to Acc-SVRG-G in \citep{zhou2021practical}, which leads to a similar rate. We also find that our correction can be used in Varag \citep{VARAG} and ANITA \citep{Anita} in the general convex case. Since the previous works  on asynchronous lock-free optimization mainly focus on the strongly convex case, we omit the discussion here.

\begin{algorithm*}[t]
	\caption{AS-Acc-SVRG: Asynchronous Sparse Accelerated SVRG}
	\label{alg:AS-Acc-SVRG}
	\renewcommand{\algorithmicrequire}{\textbf{Input:}}
	\renewcommand{\algorithmicensure}{\textbf{Initialize:}}
	\begin{algorithmic}[1]
		\Require initial guess $x_0 \in \R^d$, constant $\omega > 1$ which controls the restart frequency.
		\Ensure set $m, \vartheta, \varphi, \eta,  S, R$ according to Theorem \ref{thm:AS-Acc-SVRG}, initialize shared variable $z$ and diagonal matrix $D$.
		\For{$r=0, \ldots, R-1$}\Comment{performing restarts}
		\State $\tilde{x}_0 = z = x_r$.
		\For{$s=0, \ldots, S-1$}
		\State Compute in parallel and store $\pf{\tilde{x}_s}$.
		\While{number of samples $\leq m$} \textbf{in parallel}
		\State Sample $i$ uniformly in $[n]$ and let $T_i$ be the support of $f_i$.
		\State $\suppi{\hat{z}}=$ inconsistent read of $z$ on $T_i$.
		\State $\suppi{\hat{y}} = \vartheta\cdot \suppi{\hat{z}} + \left(1 - \vartheta\right)\cdot \suppi{\tilde{x}_s} - \varphi\cdot \suppi{D\pf{\tilde{x}_s}}.$ 
		\State $\suppi{u} = - \eta \cdot\big(\pfi{\suppi{\hat{y}}} - \pfi{\suppi{\tilde{x}_s}} + D_i\pf{\tilde{x}_s}\big)$.
		\For{$v\in T_i$}
		\State $[z]_v = [z]_v + [u]_v$. \Comment{coordinate-wise atomic write}
		\EndFor
		\EndWhile
		\State $\tilde{x}_{s+1} = \hat{y}_t$, where $\hat{y}_t$ is chosen uniformly at random among the inconsistent $\hat{y}$ in the previous epoch.
		\EndFor
		\State $x_{r+1} = \frac{1}{S}\sum_{s =0}^{S-1}{\tilde{x}_{s+1}}$.
		\EndFor
		\renewcommand{\algorithmicensure}{\textbf{Output:}}
		\Ensure $x_R$.
	\end{algorithmic}
\end{algorithm*}

\section{Asynchronous Sparse Accelerated SVRG}
\label{sec:async}
We then extend Algorithm \ref{alg:SS-Acc-SVRG} into the asynchronous setting (Algorithm \ref{alg:AS-Acc-SVRG}), and analyze it under the perturbed iterate framework proposed in \citep{ManiaPPRRJ17}. Note that Algorithm \ref{alg:AS-Acc-SVRG} degenerates into Algorithm \ref{alg:SS-Acc-SVRG} if there is only one thread (or worker).

\paragraph{Perturbed iterate analysis} let us denote the $k$th update as $\mathcal{G}_{\hat{y}_k} = \pfik{\hat{y}_k} - \pfik{\tilde{x}_s} + D_{i_k} \pf{\tilde{x}_s}$. The precise ordering of the parallel updates will be defined in the next paragraph. \citet{ManiaPPRRJ17} proposed to analyze the following virtual iterates:
\begin{equation}\label{virtual_iterate}
	z_{k+1} = z_k - \eta\cdot\mathcal{G}_{\hat{y}_k}, \text{ for } k = 0, \ldots, m-1.
\end{equation}
They are called the virtual iterates because that except for $z_0$ and $z_m$, other iterates may not exist in the shared memory due to the lock-free design. Then, the inconsistent read $\hat{z}_k$ is interpreted as a perturbed version of $z_k$, which will be formalized shortly. Note that $z_m$ is precisely $z$ in the shared memory after all the one-epoch updates are completed due to the atomic write requirement, which is critical in our analysis. 

\paragraph{Ordering the iterates} An important issue of the analysis in asynchrony is how to correctly order the updates that happen in parallel. \citet{Hogwild} increases the counter $k$ after each successful write of the update to the shared memory, and this ordering has been used in many follow-up works. Note that under this ordering, the iterates $z_k$ at \eqref{virtual_iterate} exist in the shared memory. However, this ordering is incompatible with the unbiasedness assumption (Assumption \ref{assum:indp} below), that is, enforcing the unbiasedness would require some additional overly strong assumption. \citet{ManiaPPRRJ17} addressed this issue by increasing the counter $k$ just before each inconsistent read of $z$. In this case, the unbiasedness can be simply enforced by reading all the coordinates of $z$ (not just those on the support). Although this is expensive and is not used in their implementation, the unbiasedness is enforceable under this ordering, which is thus more reasonable.  \citet{ASAGA_JMLR} further refined and simplified their analysis by proposing to increase $k$ after each inconsistent read of $z$ is completed. This modification removes the dependence of $\hat{z}_k$ on ``future'' updates, i.e., on $i_r$ for $r > k$, which significantly simplifies the analysis and leads to better speed-up conditions. See \citep{ASAGA_JMLR} for more detailed discussion on this issue. We follow the ordering of \citep{ASAGA_JMLR} to analyze Algorithm \ref{alg:AS-Acc-SVRG}. Given this ordering, the value of $\hat{z}_k$ can be explicitly described as
\[
[\hat{z}_k]_v = [z_0]_v - \eta\sum_{\substack{j\in\{0,\ldots, k-1\}\\ \text{s.t. coordinate $v$ was}\\ \text{written for $j$ before $k$}}} {[\mathcal{G}_{\hat{y}_j}]_v}.
\]
Since $\hat{y}$ is basically composing $\hat{z}$ with constant vectors, it can also be ordered as \begin{equation} \label{def_inconsistency_y}\hat{y}_k = \vartheta\hat{z}_k + \left(1 - \vartheta\right)\tilde{x}_s - \varphi D\pf{\tilde{x}_s}.
\end{equation}

\paragraph{Assumptions} The analysis in this line of work crucially relies on the following two assumptions.
\begin{assumptions}[unbiasedness]\label{assum:indp}
	$\hat{z}_k$ is independent of the sample $i_k$. Thus,  we have $\E{\mathcal{G}_{\hat{y}_k} | \hat{z}_k} = \pf{\hat{y}_k}$.
\end{assumptions}
As we mentioned before, the unbiasedness can be enforced by reading all the coordinates of $z$ while in practice one would only read those necessary coordinates. This inconsistency exists in all the follow-up works that use the revised ordering \citep{ManiaPPRRJ17,ASAGA_JMLR, ProxASAGA,MiG,JoulaniGS19, gu2020unified}, and it is currently unknown how to avoid. Another inconsistency in Algorithm \ref{alg:AS-Acc-SVRG} is that in the implementation, when a $\hat{y}_k$ is selected as the next snapshot $\tilde{x}_{s+1}$, all the coordinates of $\hat{z}_k$ are loaded. This makes the choice of $\tilde{x}_{s+1}$ not necessarily uniformly random. KroMagnon (the improved version in \citep{ASAGA_JMLR}) also has this issue. In practice, no noticeable negative impact is observed for these two inconsistencies.
\begin{assumptions}[bounded overlaps]
	\label{assum:overlap}
	There exists a uniform bound $\tau$ on the maximum number of iterations that can overlap together.
\end{assumptions}
This is a common assumption in the analysis with stale gradients. Under this assumption, the explicit effect of asynchrony can be modeled as:
\begin{equation}\label{def_inconsistency}
	\hat{z}_k = z_k + \eta\sum_{j=(k-\tau)_{+}}^{k-1}{J^k_j \mathcal{G}_{\hat{y}_j}},
\end{equation}
where $J^k_j$ is a diagonal matrix with its elements in $\{0, 1\}$. The $1$ elements indicate that $\hat{z}_k$ lacks some ``past" updates on those coordinates.

Defining the same sparsity measure as in \citep{Hogwild}: $\Delta = \frac{1}{n}\cdot\max_{v\in[d]}{|\{i: v\in T_i\}|}$, we are now ready to establish the convergence result of Algorithm \ref{alg:AS-Acc-SVRG}. We first present the following guarantee for any $\tau$ given that some upper estimation of $\tau$ is available. The proof is provided in Appendix \ref{app:thm_async}. 

\begin{theorems}
	\label{thm:AS-Acc-SVRG}
	For some $\widetilde{\tau} \geq \tau$ and any constant $\omega > 1$, we choose $m=\Theta(n), \vartheta=\frac{\sqrt{m}}{\sqrt{\kappa(1+2\sqrt{\Delta}\widetilde{\tau})} + \sqrt{m}}, \varphi = \frac{1 - \vartheta}{L}$, $\eta =  \frac{(1 - \vartheta)}{L\vartheta(1 + 2\sqrt{\Delta}\widetilde{\tau})}$ and $S= \left\lceil 2\omega\sqrt{\frac{\kappa}{m}(1+2\sqrt{\Delta}\widetilde{\tau})} \right\rceil$. For any accuracy \mbox{$\epsilon > 0$,} we restart $R= O\big(\log{\frac{f(x_0) - f(\xs)}{\epsilon}}\big)$ rounds. In this case, Algorithm \ref{alg:AS-Acc-SVRG}   outputs $x_R$ satisfying $\E{f(x_R)} - f(\xs) \leq \epsilon$ in oracle complexity
	\[
	O\left( \max{\left\{n, \sqrt{\kappa n(1+2\sqrt{\Delta}\widetilde{\tau})} \right\}}\log{\frac{f(x_0) - f(\xs)}{\epsilon}}\right).
	\]
\end{theorems}

Based on this theorem, it is direct to identify the region of $\tau$ where a theoretical linear speed-up is achievable. The proof is straightforward and is thus omitted.

\begin{corollarys}[Speed-up condition]
	\label{cor:speed-up}
	In Theorem~\ref{thm:AS-Acc-SVRG}, suppose $\tau \leq O\left(\frac{1}{\sqrt{\Delta}}\max{\left\lbrace\frac{n}{\kappa}, 1\right\rbrace}\right)$. Then, setting $\widetilde{\tau} = O\left(\frac{1}{\sqrt{\Delta}}\max{\left\lbrace\frac{n}{\kappa}, 1\right\rbrace}\right)$, we have $S = O\left(\max{\left\lbrace 1, \sqrt{\frac{\kappa}{n}}\right\rbrace}\right)$, which leads to the total oracle complexity $\text{\#grad} = R\cdot S\cdot(n+2m) = O\big(\max{\left\{n, \sqrt{\kappa n}\right\}}\log{\frac{f(x_0) - f(\xs)}{\epsilon}}\big)$.
\end{corollarys}
%\begin{proof}
	%	Note that $\widetilde{\tau} = O\big(\frac{1}{\sqrt{\Delta}}\max{\left\lbrace\frac{n}{\kappa}, 1\right\rbrace}\big)$ implies that $1 + 2\sqrt{\Delta}\widetilde{\tau} = O\left(\max{\left\lbrace\frac{n}{\kappa}, 1\right\rbrace}\right)$. In this case, it holds that $S = O\left(\max{\left\lbrace 1, \sqrt{\frac{\kappa}{n}}\right\rbrace}\right)$, and then the total oracle complexity of Algorithm \ref{alg:AS-Acc-SVRG} is 
	%	$
	%	\text{\#grad} = R\cdot S\cdot(n+2m) = O\big(\max{\left\{n, \sqrt{\kappa n}\right\}}\log{\frac{f(x_0) - f(\xs)}{\epsilon}}\big).
	%	$
	%\end{proof}

Note that the construction of Algorithm \ref{alg:AS-Acc-SVRG} naturally enforces that $\tau \leq m$. Hence, the precise linear speed-up condition of AS-Acc-SVRG is that $\tau = O(n)$ and $\tau = O(\frac{1}{\sqrt{\Delta}}\max{\left\lbrace\frac{n}{\kappa}, 1\right\rbrace})$, which is identical to that of ASAGA (cf., Corollary 9 in \citep{ASAGA_JMLR}) and slightly better than that of KroMagnon (cf., Corollary~18 in \citep{ASAGA_JMLR}).

\subsection{Some Insights about the Asynchronous Acceleration}
Let us first consider the serial case (Algorithm \ref{alg:SS-Acc-SVRG}). Observe that in one epoch, the iterate $y_k$ is basically composing $z_k$ with constant vectors, we can equivalently write the update (Step \ref{alg-step:mirror_c}) as $y_{k+1} = y_k - \eta\vartheta \cdot \mathcal{G}_{y_k}$. Thus, the inner loop of Algorithm \ref{alg:SS-Acc-SVRG} is identical to that of (sparse) SVRG. The difference is that at the end of each epoch, when the snapshot $\tilde{x}$ is changed, an offset (or momentum) is added to the iterate. This has been similarly observed for accelerated SVRG in \citep{SSNM}. Note that since the sequence $z$ appears in the potential function, the current formulation of Algorithm \ref{alg:SS-Acc-SVRG} allows a cleaner analysis. Then, in the asynchronous case, the inner loop of Algorithm~\ref{alg:AS-Acc-SVRG} can also be equivalently written as the updates of asynchronous SVRG, and the momentum is added at the end of each epoch. This gives us some insights about the identical speed-up condition in Corollary \ref{cor:speed-up}: Since the asynchronous perturbation only affects the inner loop, the momentum is almost uncorrupted, unlike the cases of noisy gradient oracle. That is, the asynchrony only corrupts the ``non-accelerated part'' of Algorithm~\ref{alg:AS-Acc-SVRG}, which thus leads to the same speed-up condition as the non-accelerated methods.

\bgroup
\def\arraystretch{1.5}
\begin{table*}[t]
	\scriptsize
	\caption{Summary of the datasets. Density is the ratio of non-zero elements.} \label{data-table}
	\begin{center}
		\begin{tabular}{|c|l|c|c|c|c|c|c|}
			\hline \textbf{Scale}  &\textbf{Dataset}  &$\boldsymbol{n}$  &$\boldsymbol{d}$
			&$\boldsymbol{\mu}$
			&\textbf{Density} &$\boldsymbol{\Delta}$  &\textbf{Description} \\
			\hhline{|========|}
			\multirow{4}{*}{Small} &Synthetic &\SI{100000}{} &\SI{100000}{} &$10^{-7}$  &$10^{-5}$ &$10^{-5}																																									$ &Identity data matrix with random labels\\
			&\textsf{KDD2010.small}  &\SI{70000}{} &\SI{29890095}{} &$10^{-7}$ &$10^{-6}$ &$0.15$ &The first \SI{70000}{} data samples of \textsf{KDD2010} \\
			&\textsf{RCV1.train}  &\SI{20242}{} &\SI{47236}{}&$10^{-6}$ &$1.6 \cdot 10^{-3}$ &$0.42$ & \\
			&\textsf{News20} &\SI{19996}{} &\SI{1355191}{} &$10^{-6}$ &$3.4 \cdot 10^{-4}$ &$0.93$ & \\
			\hline \multirow{3}{*}{Large} &\textsf{KDD2010} &\SI{19264097}{} &\SI{29890095}{} &$10^{-10}$ &$10^{-6}$ &$0.16$& \\
			&\textsf{RCV1.full}  &\SI{697641}{} &\SI{47236}{} &$10^{-9}$ &$1.5 \cdot 10^{-3}$ &$0.43$ &Combined test and train sets of \textsf{RCV1} \\
			&\textsf{Avazu-site.train}  &\SI{23567843}{} &\SI{999962}{} &$10^{-10}$ &$1.5 \cdot 10^{-5}$ &$0.96$ & \\
			\hline
		\end{tabular}
	\end{center}
\end{table*}
\egroup

\section{Experiments}
We present numerical results for the proposed scheme on optimizing the $\ell_2$-logistic regression problem:
\begin{equation} \label{exp-obj}
	f(x) = \frac{1}{n} \sum_{i=1}^n {\log{\big(1 + \exp{(-b_i\innr{a_i,x})}\big)}} + \frac{\mu}{2}\norm{x}^2,
\end{equation}
where $a_i \in \R^d$, $b_i \in \lbrace -1, +1\rbrace$, $i\in[n]$ are the data samples and $\mu$ is the regularization parameter. We use the datasets from LIBSVM website \citep{LIBSVM}, including \textsf{KDD2010} \citep{yu2010feature}, \textsf{RCV1} \citep{lewis2004rcv1}, \textsf{News20} \citep{keerthi2005modified}, \textsf{Avazu} \citep{JuanZCL16}. All the datasets are normalized to ensure a precise control on $\kappa$. We focus on the ill-conditioned case where $\kappa \gg n$ (the case where acceleration is effective). The dataset descriptions and the choices of $\mu$ are provided in Table \ref{data-table}. We conduct serial experiments on the small datasets and asynchronous experiments on the large ones. The asynchronous experiments were conducted on a multi-core HPC. Detailed setup can be found in Appendix~\ref{app:exp}.

Before presenting the empirical results, let us discuss two subtleties in the implementation:
\begin{itemize}[leftmargin=14pt]
	\item If each $f_i(x) = \log{\big(1 + \exp{(-b_i\innr{a_i,x})}\big)} + \frac{\mu}{2}\norm{x}^2$, then it is supported on every coordinate due to the $\ell_2$-regularization. Following \citep{ASAGA_JMLR}, we sparsify the gradient of the regularization term as $\mu D_i x$; or equivalently, $f_i(x)=\log{\big(1 + \exp{(-b_i\innr{a_i,x})}\big)} + \frac{\mu}{2}\innr{x, D_ix}$. Clearly, this $f_i$ also sums up to the objective \eqref{exp-obj}. The difference is that the Lipschitz constant of $f_i$ will be larger, i.e., from $0.25 + \mu$ to at most $0.25 + \mu n$. Since we focus on the ill-conditioned case ($L \gg \mu n$), this modification will not have significant effect.
	
	\item In Theorems \ref{thm:SS-Acc-SVRG} and \ref{thm:AS-Acc-SVRG}, all the parameters have been chosen optimally in our analysis except for $\omega$, which controls the restart frequency. Despite all the theoretical benefits of the restart framework mentioned in Section \ref{sec:remarks}, in practice, we do not find performing restarts leads to a faster convergence. Thus, we do not perform restarts in our experiments (or equivalently, we choose a relatively large $\omega = 50$). Clearly, this choice will not make our method a ``heuristic'' since the theorems hold for any $\omega > 1$. Detailed discussion is given in Appendix \ref{app:restart}.
\end{itemize}

Additional experiments for verifying the $\sqrt{\kappa}$ dependence and a sanity check for our implementation is included in Appendix \ref{app:justify} and \ref{app:sanity}, respectively.

\subsection{The Effectiveness of Sparse Variance Correction}

\begin{figure*}[t]
	\begin{center}
		\begin{subfigure}{0.245\linewidth}
			\includegraphics[width=\linewidth]{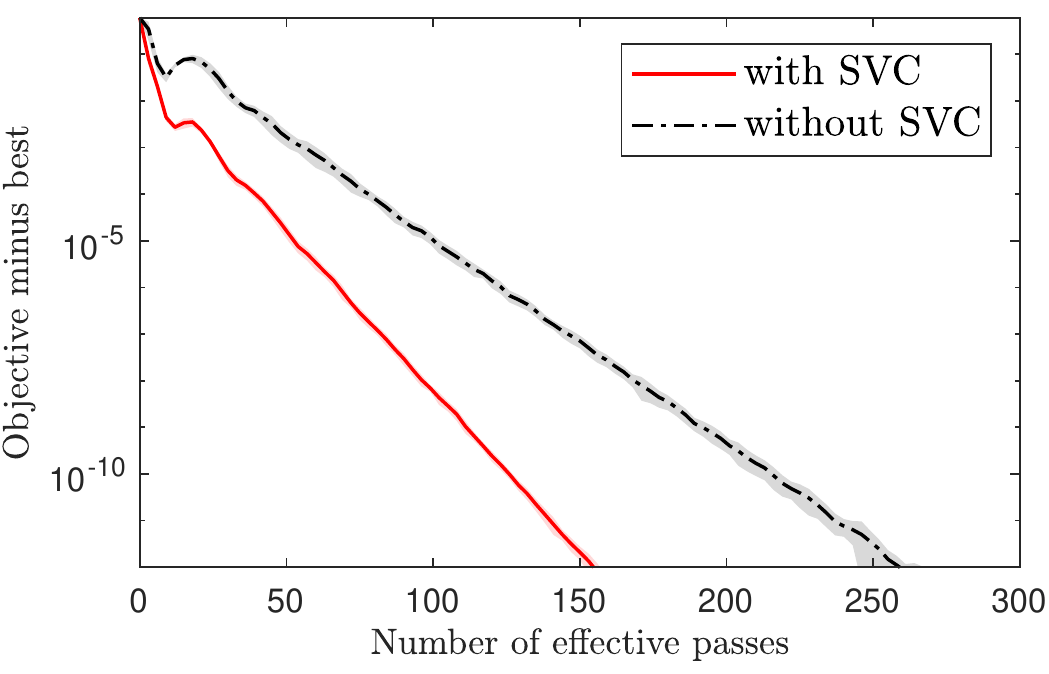}
			\caption{Synthetic}
			\label{fig:SVC_synthetic}
		\end{subfigure}
		\begin{subfigure}{0.245\linewidth}
			\includegraphics[width=\linewidth]{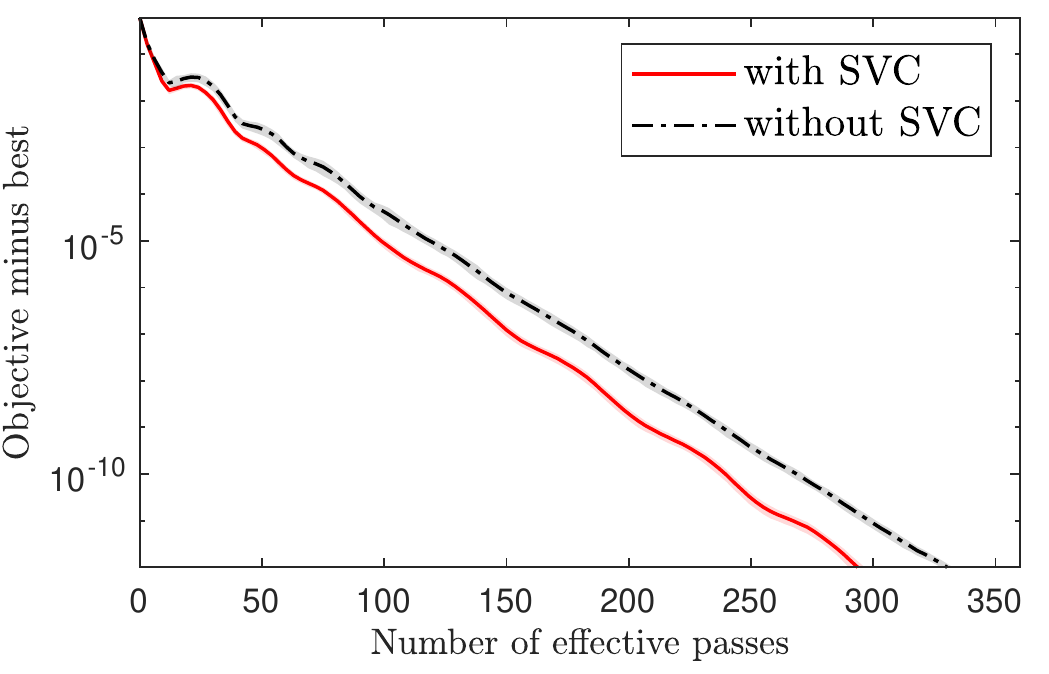}
			\caption{\textsf{KDD2010.small}}
			\label{fig:SVC_kdd10_small}
		\end{subfigure}
		\begin{subfigure}{0.245\linewidth}
			\includegraphics[width=\linewidth]{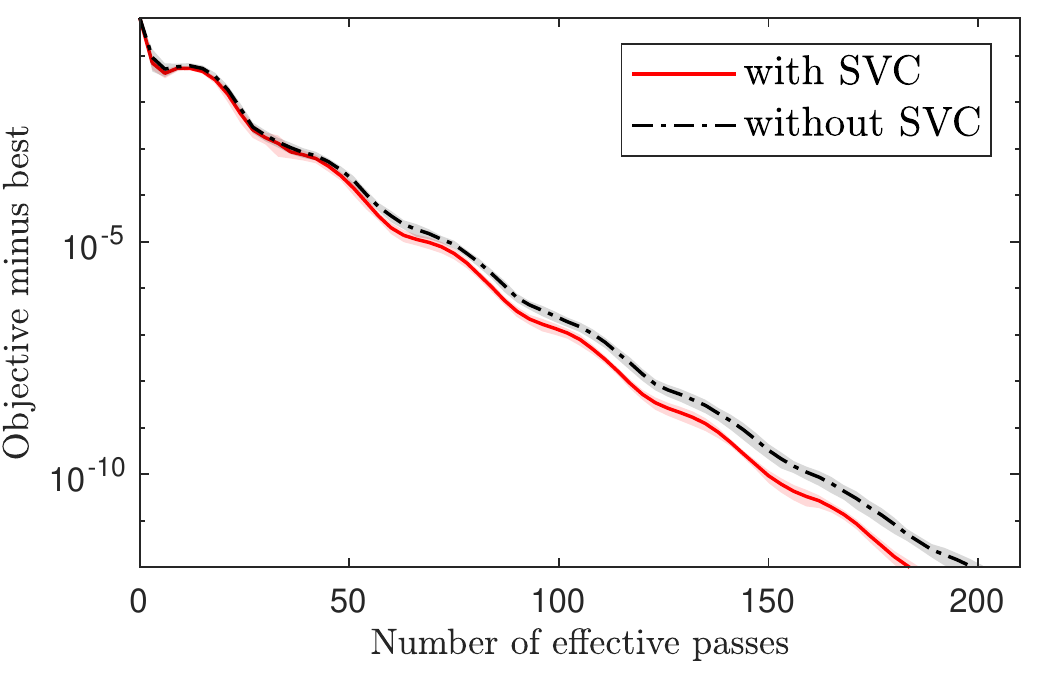}
			\caption{\textsf{RCV1.train}}
			\label{fig:SVC_rcv1_train}
		\end{subfigure}
		\begin{subfigure}{0.245\linewidth}
			\includegraphics[width=\linewidth]{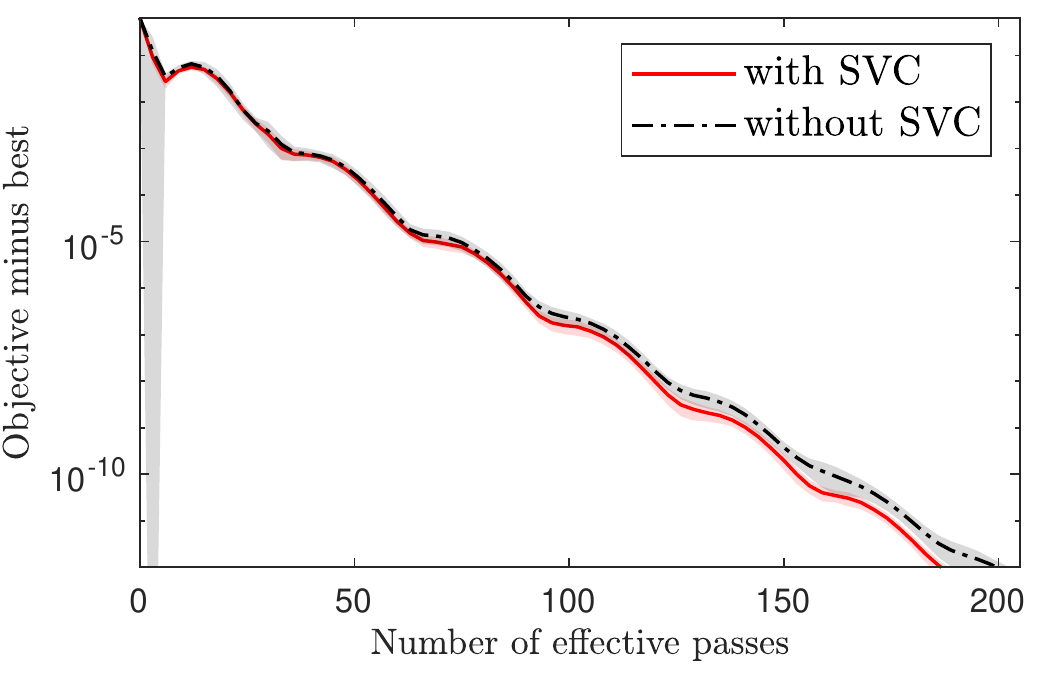}
			\caption{\textsf{News20}}
			\label{fig:SVC_news20}
		\end{subfigure}
		\caption{Ablation study for the practical effect of sparse variance correction (abbreviated as SVC in the legends). Run 10 seeds. Shaded bands indicate $\pm 1$ standard deviation.}
		\label{fig:SVC}
	\end{center}
\end{figure*}

We first study the practical effect of the correction in the serial case. In Figure \ref{fig:SVC}, ``with SVC'' is Algorithm~\ref{alg:SS-Acc-SVRG} and ``without SVC'' is a naive extension of the dense version of Algorithm~\ref{alg:SS-Acc-SVRG} into the sparse case (i.e., simply using a sparse estimator \eqref{S-SVRG-estimator}), which suffers from the same sparsity issue as in \citep{MiG}. For the two variants, we chose the same parameters specified in Theorem \ref{alg:SS-Acc-SVRG} to conduct an ablation study. In theory, we would expect the sparsity issue to be more severe if $D$ is closer to $n I_d$ (i.e., more sparse). Note that by definition, $\Delta = \max_{i\in[d]} {D_{ii}^{-1}}$, and thus $D \succeq \Delta^{-1} I_d$. Hence, smaller $\Delta$ indicates that $D$ is closer to $n I_d$. From Figure \ref{fig:SVC_news20} to \ref{fig:SVC_synthetic}, $\Delta$ is decreasing and we observe more improvement from the correction. The improvement is consistent in our experiments, which justifies the effectiveness of sparse variance correction.

\subsection{Sparse Estimator v.s. Lagged Update}
We then examine the running time improvement from adopting the sparse estimator compared with using the lagged update technique. Lagged update technique handles sparsity by maintaining a last seen iteration for each coordinate. When a coordinate is involved in the current iteration, it computes an accumulated update from the last seen iteration in closed form. Such computation is dropped when adopting a sparse estimator, which is the source of running time improvement. Moreover, it is extremely difficult to extend the lagged update technique into the asynchronous setting as discussed in Appendix E in \citep{ASAGA_JMLR}. In Figure \ref{fig:LU}, we compare SS-Acc-SVRG (Algorithm \ref{alg:SS-Acc-SVRG}) with lagged update implementations of (dense) SS-Acc-SVRG and Katyusha. Their default parameters were used. Note that the lagged update technique is much trickier to implement, especially for Katyusha. We need to hand compute the closed-form solution of some complicated constant recursive sequence for the accumulated update. Plots with respect to effective passes are provided in Appendix \ref{app:justify}, in which SS-Acc-SVRG and Katyusha show similar performance.

\begin{figure*}[t]
	\begin{center}
		\begin{subfigure}{0.245\linewidth}
			\includegraphics[width=\linewidth]{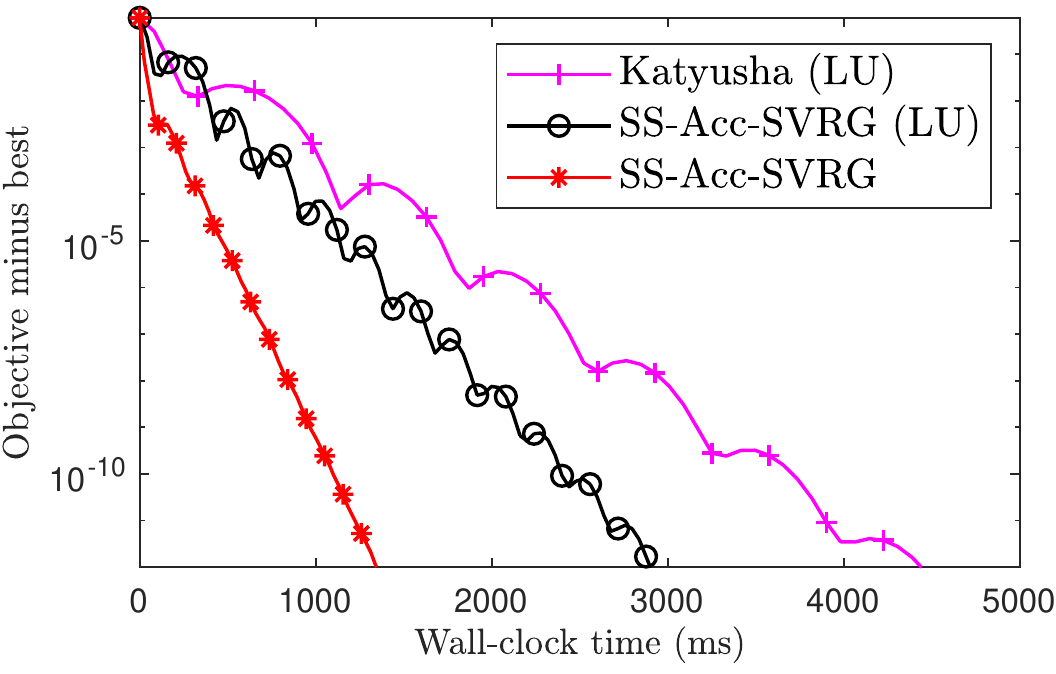}
			\caption{Synthetic}
			\label{fig:LU_synthetic}
		\end{subfigure}
		\begin{subfigure}{0.245\linewidth}
			\includegraphics[width=\linewidth]{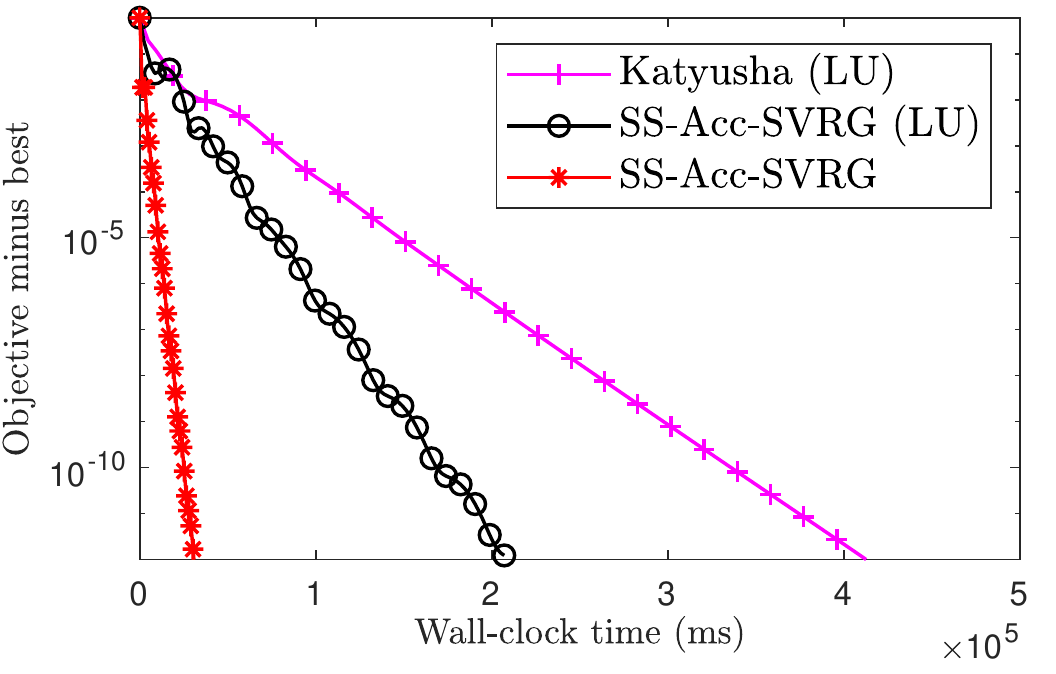}
			\caption{\textsf{KDD2010.small}}
			\label{fig:LU_kdd10_small}
		\end{subfigure}
		\begin{subfigure}{0.245\linewidth}
			\includegraphics[width=\linewidth]{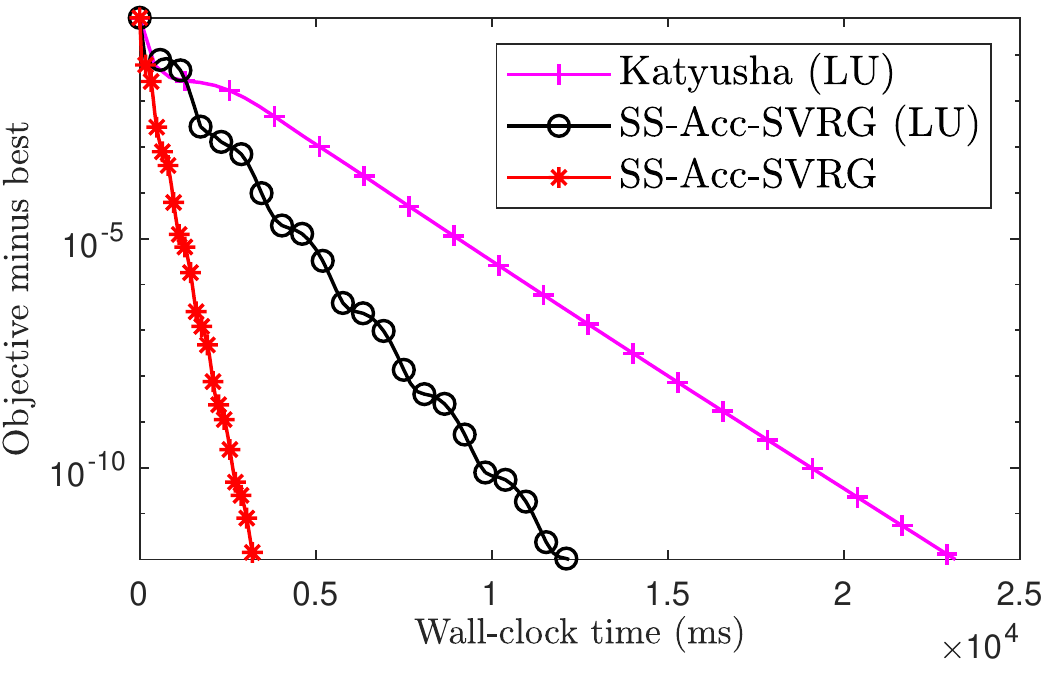}
			\caption{\textsf{RCV1.train}}
			\label{fig:LU_rcv1_train}
		\end{subfigure}
		\begin{subfigure}{0.245\linewidth}
			\includegraphics[width=\linewidth]{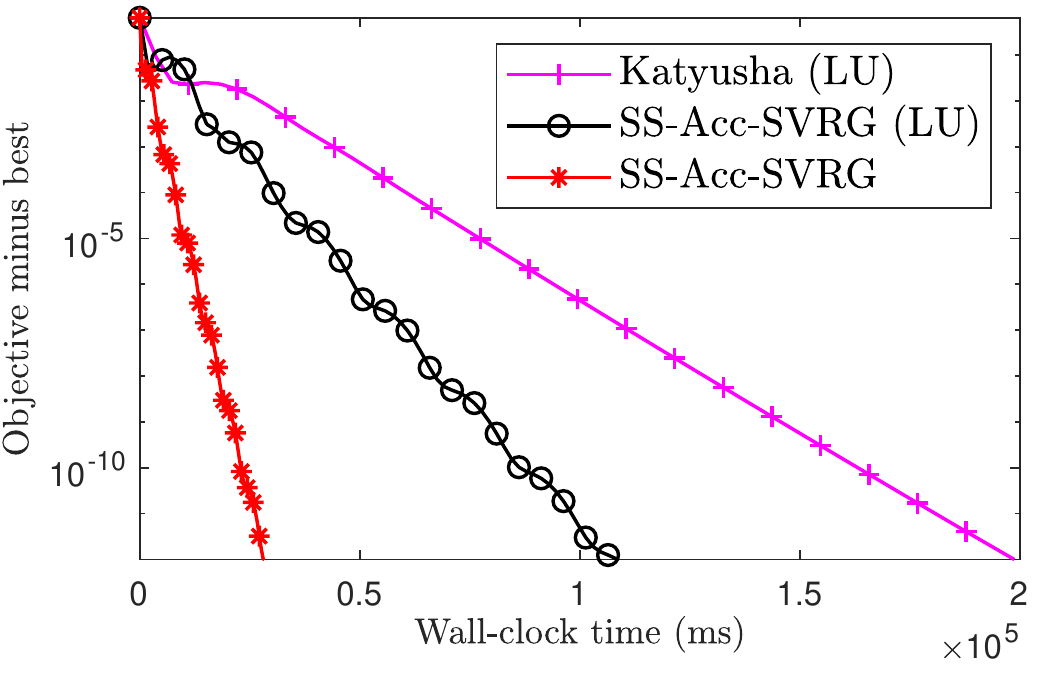}
			\caption{\textsf{News20}}
			\label{fig:LU_news20}
		\end{subfigure}
		\caption{Running time comparison between using sparse gradient estimator and lagged update (abbreviated as LU in the legends). The wall-clock time and objective value are averaged over 10 runs. }
		\label{fig:LU}
	\end{center}
\end{figure*}
\begin{figure*}[t]
	\begin{center}
		\begin{subfigure}{0.27\linewidth}
			\includegraphics[width=\linewidth]{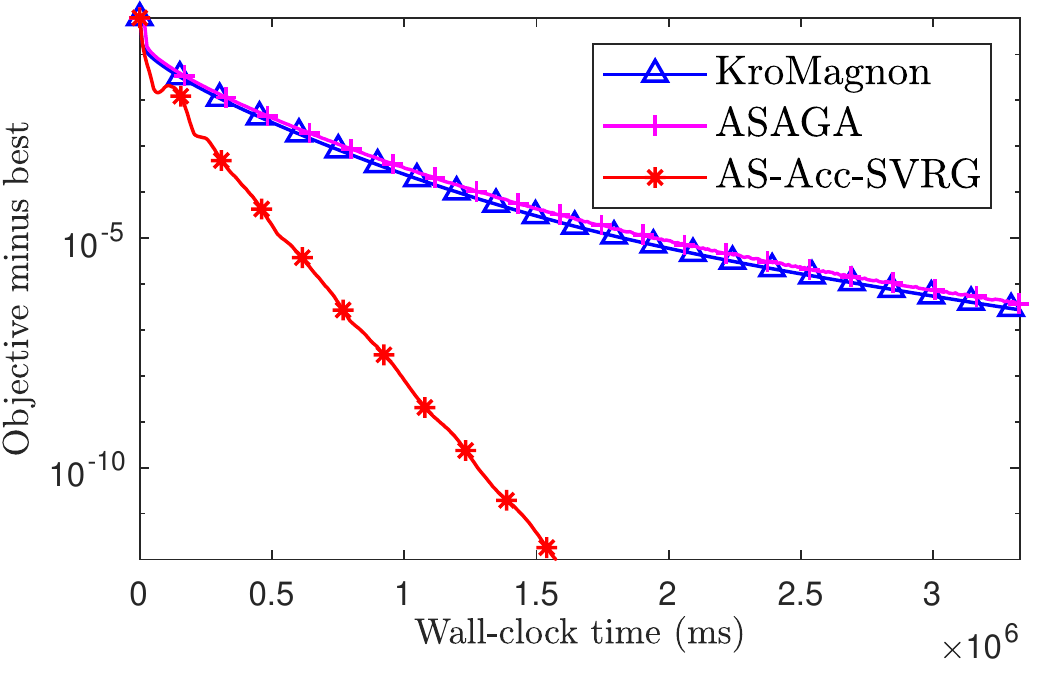}
			\caption{\textsf{KDD2010}, $20$ threads}
			\label{fig:AS_kdd10}
		\end{subfigure}
		\quad\quad
		\begin{subfigure}{0.27\linewidth}
			\includegraphics[width=\linewidth]{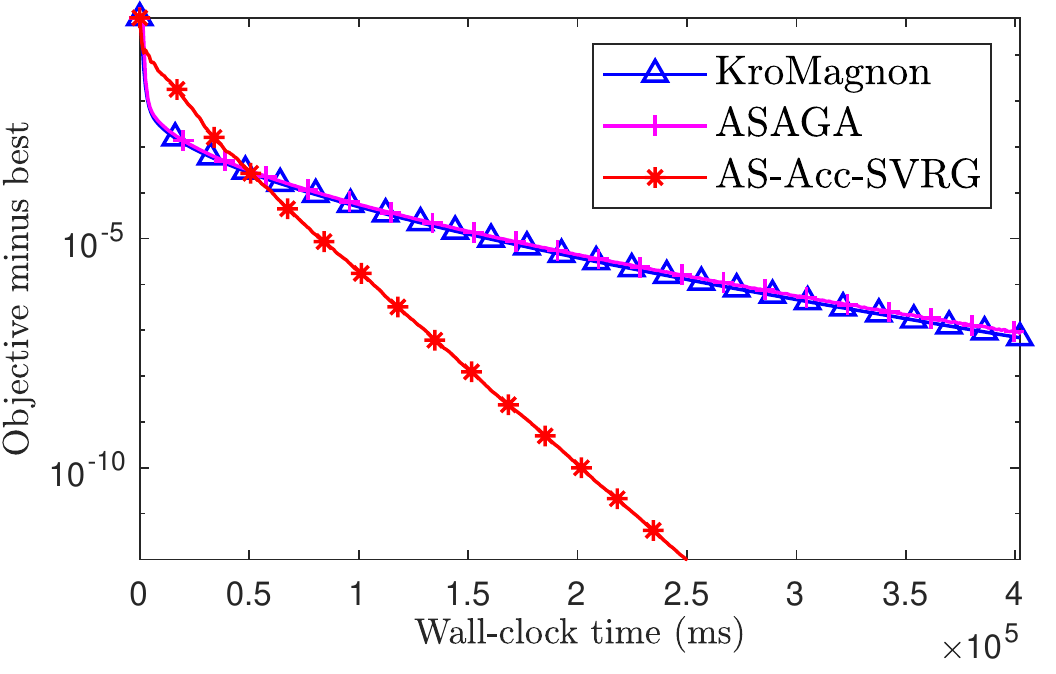}
			\caption{\textsf{RCV1.full}, $20$ threads}
			\label{fig:AS_rcv1}
		\end{subfigure}
		\quad\quad
		\begin{subfigure}{0.27\linewidth}
			\includegraphics[width=\linewidth]{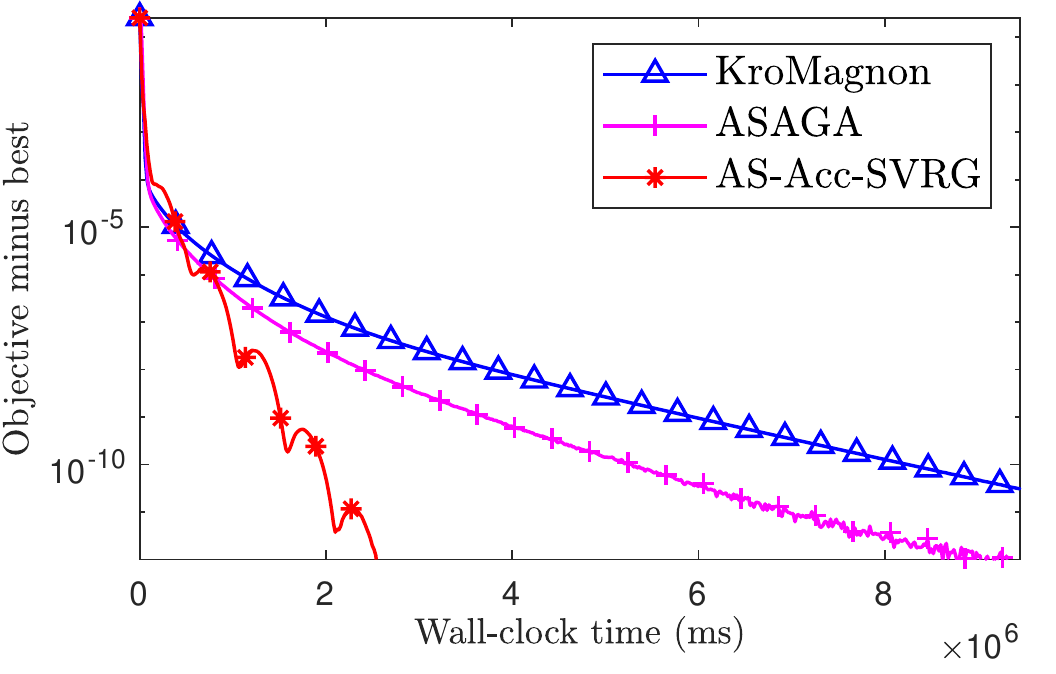}
			\caption{\textsf{Avazu-site.train}, $20$ threads}
			\label{fig:AS_avazu}
		\end{subfigure} \\
		\,\,
		\begin{subfigure}{0.27\linewidth}
			\includegraphics[width=0.9\linewidth]{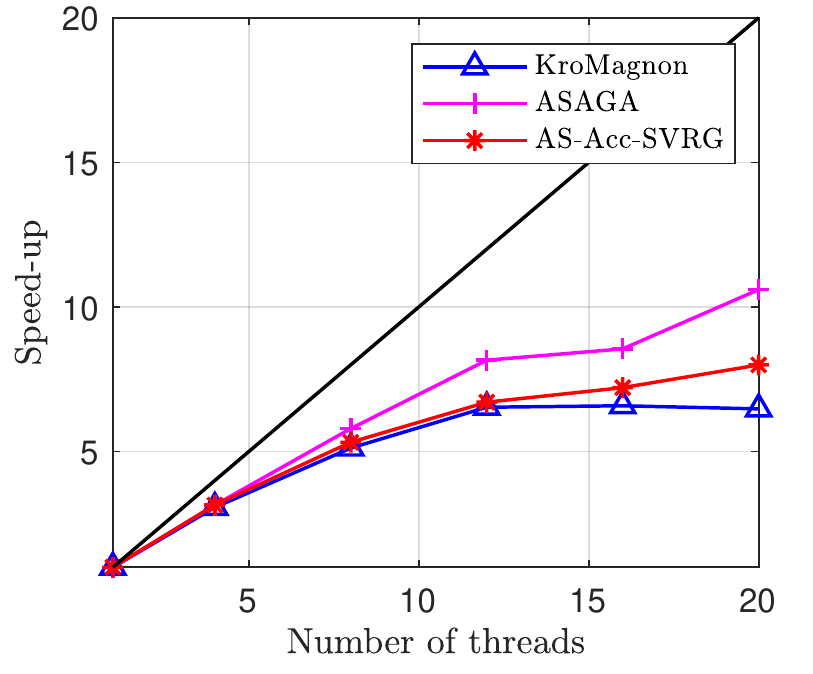}
			\centering
			\caption{\textsf{KDD2010}, speed-up}
			\label{fig:SU_kdd10}
		\end{subfigure}\quad\quad\,
		\begin{subfigure}{0.27\linewidth}
			\includegraphics[width=0.9\linewidth]{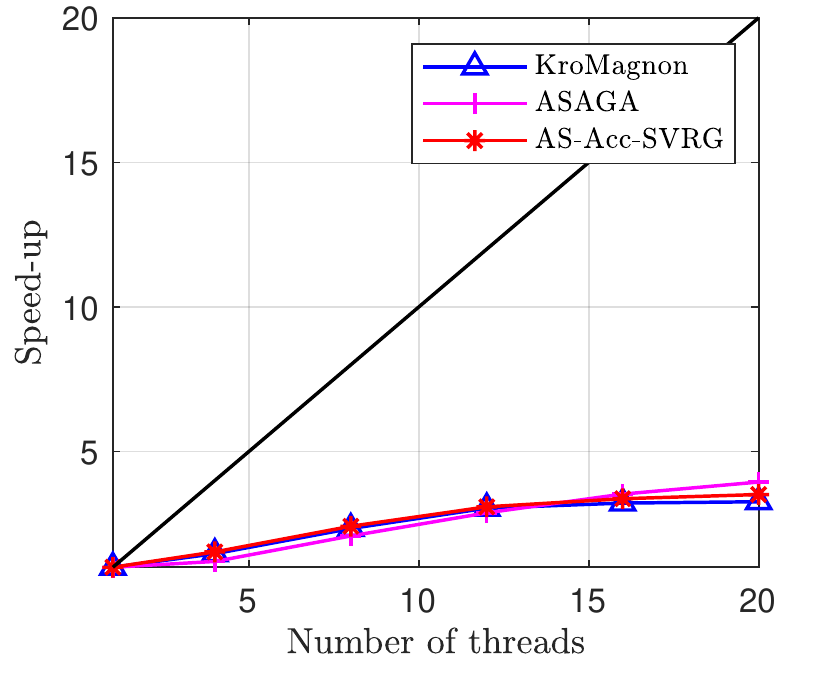}
			\centering
			\caption{\textsf{RCV1.full}, speed-up}
			\label{fig:SU_rcv1}
		\end{subfigure}\quad\quad\,
		\begin{subfigure}{0.27\linewidth}
			\includegraphics[width=0.9\linewidth]{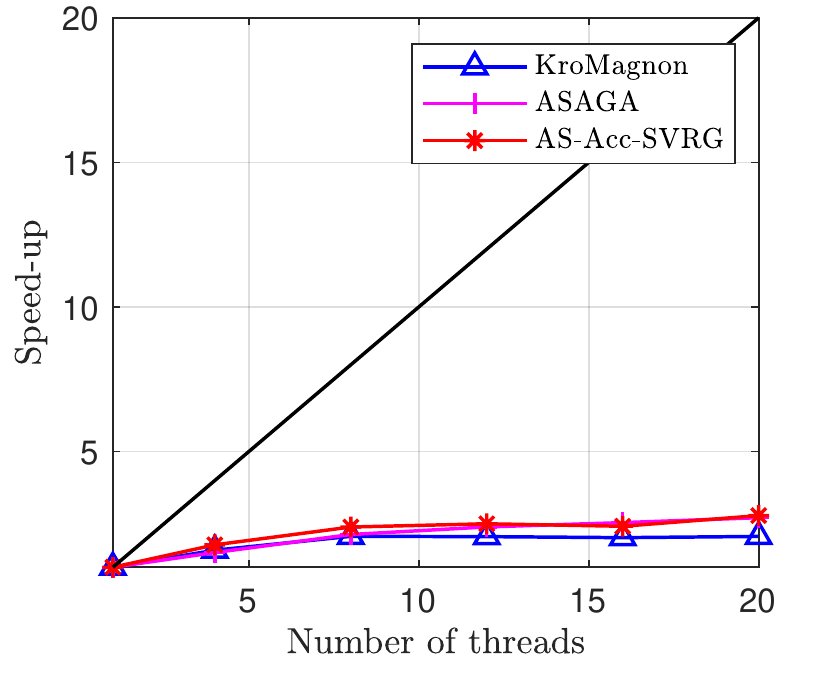}
			\centering
			\caption{\textsf{Avazu-site.train}, speed-up}
			\label{fig:SU_avazu}
		\end{subfigure}
		\caption{Convergence and speed-up for asynchronous sparse methods. Speed-up is the improvement on the wall-clock time to achieve $10^{-5}$ sub-optimality relative to using a single thread.}
		\label{fig:AS_SU}
	\end{center}
\end{figure*}

\subsection{Asynchronous Experiments}
We compare the practical convergence and speed-up of AS-Acc-SVRG (Algorithm \ref{alg:AS-Acc-SVRG}) with KroMagnon and ASAGA in Figure \ref{fig:AS_SU}.  We do not compare with the empirical method MiG in \citep{MiG}, which requires us to tune two highly correlated parameters with only limited insights. This is expensive or even prohibited for large scale tasks. For the compared methods, we only tune the $\tau$-related constants in their theoretical parameter settings. That is, we tune the constant $1 + 2\sqrt{\Delta}\tilde{\tau}$ in Theorem \ref{thm:AS-Acc-SVRG} for AS-Acc-SVRG and the constant $c$ in the step size $\frac{1}{cL}$ for KroMagnon and ASAGA. In fact, in all the experiments, we simply fixed the constant to $1$ for AS-Acc-SVRG (the same parameters as the serial case), which worked smoothly. The main tuning effort was devoted to KroMagnon and ASAGA, and we tried to choose their step sizes as large as possible. The  detailed choices can be found in Appendix \ref{app:exp}.  Due to the scale of the problems, we only conducted a single run. From Figure~\ref{fig:AS_SU}, we see significant improvement of AS-Acc-SVRG for ill-conditioned tasks and similar practical speed-ups among the three methods, which verifies Theorem \ref{thm:AS-Acc-SVRG} and Corollary \ref{cor:speed-up}. We also observe a strong correlation between the practical speed-up and $\Delta$, which is predicted by the theoretical $O(1/\sqrt{\Delta})$ dependence. It seems that we cannot reproduce the speed-up results in \citep{ASAGA_JMLR} on \textsf{RCV1.full} dataset. We believe that it is due to the usage of different programming languages, that we used a low-level heavily-optimized language \CC\,and they used high-level language Scalar. \citet{ProxASAGA} also used \CC\,implementation and we observe a similar $10\times$ speed-up of ASAGA on \textsf{KDD2010} dataset.

\section{Conclusion}
In this work, we proposed a new asynchronous accelerated SVRG method which achieves the optimal oracle complexity under the perturbed iterate framework. We show that it requires the same linear speed-up condition as the non-accelerated methods. Empirical results justified our findings.

A limitation of this work is that it does not support proximal operators, and we are currently not sure how to leverage the techniques in \citep{ProxASAGA} into our method. The main difficulty is on the choice of potential function. 

Future directions could be (i) a more fine-grained analysis on $\Delta$ since currently $\Delta$ is a constant in all the datasets (except for the synthetic one), and (ii) designing asynchronous accelerated SAGA based on \citep{SSNM} since SVRG variants still require at least one synchronization per epoch.
\bibliography{Sparse_Acc}

\begin{thebibliography}{}

\bibitem[Allen-Zhu, 2017]{allen2017katyusha}
Allen-Zhu, Z. (2017).
\newblock {Katyusha: The First Direct Acceleration of Stochastic Gradient
  Methods}.
\newblock {\em Journal of Machine Learning Research}, 18(1):8194--8244.

\bibitem[Chang and Lin, 2011]{LIBSVM}
Chang, C.-C. and Lin, C.-J. (2011).
\newblock {LIBSVM}: A library for support vector machines.
\newblock {\em ACM Transactions on Intelligent Systems and Technology},
  2:27:1--27:27.
\newblock Software available at \url{http://www.csie.ntu.edu.tw/~cjlin/libsvm}.

\bibitem[Cyrus et~al., 2018]{RMM}
Cyrus, S., Hu, B., Van~Scoy, B., and Lessard, L. (2018).
\newblock A robust accelerated optimization algorithm for strongly convex
  functions.
\newblock In {\em Annual American Control Conference (ACC)}, pages 1376--1381.
  IEEE.

\bibitem[Defazio et~al., 2014]{SAGA}
Defazio, A., Bach, F.~R., and Lacoste{-}Julien, S. (2014).
\newblock {SAGA: A Fast Incremental Gradient Method With Support for
  Non-Strongly Convex Composite Objectives}.
\newblock In {\em Advances in Neural Information Processing Systems}, pages
  1646--1654.

\bibitem[Devolder et~al., 2014]{devolder2014first}
Devolder, O., Glineur, F., and Nesterov, Y. (2014).
\newblock {First-order methods of smooth convex optimization with inexact
  oracle}.
\newblock {\em Mathematical Programming}, 146(1):37--75.

\bibitem[Dubois-Taine et~al., 2021]{dubois2021svrg}
Dubois-Taine, B., Vaswani, S., Babanezhad, R., Schmidt, M., and Lacoste-Julien,
  S. (2021).
\newblock {SVRG Meets AdaGrad: Painless Variance Reduction}.
\newblock {\em arXiv preprint arXiv:2102.09645}.

\bibitem[Duchi et~al., 2013]{DuchiJM13}
Duchi, J.~C., Jordan, M.~I., and McMahan, H.~B. (2013).
\newblock {Estimation, Optimization, and Parallelism when Data is Sparse}.
\newblock In {\em Advances in Neural Information Processing Systems}, pages
  2832--2840.

\bibitem[Fang et~al., 2018]{fang2018accelerating}
Fang, C., Huang, Y., and Lin, Z. (2018).
\newblock {Accelerating asynchronous algorithms for convex optimization by
  momentum compensation}.
\newblock {\em arXiv preprint arXiv:1802.09747}.

\bibitem[Gu et~al., 2020]{gu2020unified}
Gu, B., Xian, W., Huo, Z., Deng, C., and Huang, H. (2020).
\newblock {A Unified q-Memorization Framework for Asynchronous Stochastic
  Optimization}.
\newblock {\em Journal of Machine Learning Research}, 21:190--1.

\bibitem[Hannah et~al., 2019]{A2BCD}
Hannah, R., Feng, F., and Yin, W. (2019).
\newblock {A2{BCD}: Asynchronous Acceleration with Optimal Complexity}.
\newblock In {\em International Conference on Learning Representations}.

\bibitem[Hofmann et~al., 2015]{HofmannLLM15}
Hofmann, T., Lucchi, A., Lacoste{-}Julien, S., and McWilliams, B. (2015).
\newblock {Variance Reduced Stochastic Gradient Descent with Neighbors}.
\newblock In {\em Advances in Neural Information Processing Systems}, pages
  2305--2313.

\bibitem[Johnson and Zhang, 2013]{SVRG}
Johnson, R. and Zhang, T. (2013).
\newblock {Accelerating Stochastic Gradient Descent using Predictive Variance
  Reduction}.
\newblock In {\em Advances in Neural Information Processing Systems}, pages
  315--323.

\bibitem[Joulani et~al., 2019]{JoulaniGS19}
Joulani, P., Gy{\"{o}}rgy, A., and Szepesv{\'{a}}ri, C. (2019).
\newblock {Think out of the "Box": Generically-Constrained Asynchronous
  Composite Optimization and Hedging}.
\newblock In {\em Advances in Neural Information Processing Systems}, pages
  12225--12235.

\bibitem[Juan et~al., 2016]{JuanZCL16}
Juan, Y., Zhuang, Y., Chin, W., and Lin, C. (2016).
\newblock {Field-aware Factorization Machines for {CTR} Prediction}.
\newblock In {\em Proceedings of the 10th {ACM} Conference on Recommender
  Systems}, pages 43--50. {ACM}.

\bibitem[Karimi et~al., 2016]{karimi2016linear}
Karimi, H., Nutini, J., and Schmidt, M. (2016).
\newblock {Linear Convergence of Gradient and Proximal-Gradient Methods Under
  the Polyak-{\L}ojasiewicz Condition}.
\newblock In {\em Joint European Conference on Machine Learning and Knowledge
  Discovery in Databases}, pages 795--811. Springer.

\bibitem[Keerthi et~al., 2005]{keerthi2005modified}
Keerthi, S.~S., DeCoste, D., and Joachims, T. (2005).
\newblock {A Modified Finite Newton Method for Fast Solution of Large Scale
  Linear SVMs}.
\newblock {\em Journal of Machine Learning Research}, 6:341--361.

\bibitem[Kone{\v{c}}n{\`y} and Richt{\'a}rik, 2013]{S2GD}
Kone{\v{c}}n{\`y}, J. and Richt{\'a}rik, P. (2013).
\newblock Semi-stochastic gradient descent methods.
\newblock {\em arXiv preprint arXiv:1312.1666}.

\bibitem[Lan et~al., 2019]{VARAG}
Lan, G., Li, Z., and Zhou, Y. (2019).
\newblock A unified variance-reduced accelerated gradient method for convex
  optimization.
\newblock In {\em Advances in Neural Information Processing Systems},
  volume~32, pages 10462--10472.

\bibitem[Leblond et~al., 2017]{ASAGA}
Leblond, R., Pedregosa, F., and Lacoste{-}Julien, S. (2017).
\newblock {ASAGA: Asynchronous Parallel SAGA}.
\newblock In {\em Proceedings of the Twentieth International Conference on
  Artificial Intelligence and Statistics}, volume~54, pages 46--54.

\bibitem[Leblond et~al., 2018]{ASAGA_JMLR}
Leblond, R., Pedregosa, F., and Lacoste{-}Julien, S. (2018).
\newblock {Improved Asynchronous Parallel Optimization Analysis for Stochastic
  Incremental Methods}.
\newblock {\em Journal of Machine Learning Research}, 19:81:1--81:68.

\bibitem[Lewis et~al., 2004]{lewis2004rcv1}
Lewis, D.~D., Yang, Y., Russell-Rose, T., and Li, F. (2004).
\newblock {{RCV1:} {A} New Benchmark Collection for Text Categorization
  Research}.
\newblock {\em Journal of Machine Learning Research}, 5:361--397.

\bibitem[Li, 2021]{Anita}
Li, Z. (2021).
\newblock {ANITA: An Optimal Loopless Accelerated Variance-Reduced Gradient
  Method}.
\newblock {\em arXiv preprint arXiv:2103.11333}.

\bibitem[Lian et~al., 2015]{LianHLL15}
Lian, X., Huang, Y., Li, Y., and Liu, J. (2015).
\newblock {Asynchronous Parallel Stochastic Gradient for Nonconvex
  Optimization}.
\newblock In {\em Advances in Neural Information Processing Systems}, pages
  2737--2745.

\bibitem[Mania et~al., 2017]{ManiaPPRRJ17}
Mania, H., Pan, X., Papailiopoulos, D.~S., Recht, B., Ramchandran, K., and
  Jordan, M.~I. (2017).
\newblock {Perturbed Iterate Analysis for Asynchronous Stochastic
  Optimization}.
\newblock {\em SIAM Journal on Optimization}, 27(4):2202--2229.

\bibitem[Nesterov, 1983]{AGD1}
Nesterov, Y. (1983).
\newblock {A method for solving the convex programming problem with convergence
  rate $O(1/k^2)$}.
\newblock In {\em Dokl. akad. nauk Sssr}, volume 269, pages 543--547.

\bibitem[Nesterov, 2003]{AGD2}
Nesterov, Y. (2003).
\newblock {\em Introductory lectures on convex optimization: A basic course},
  volume~87.
\newblock Springer Science \& Business Media.

\bibitem[Nesterov, 2005]{nesterov2005smooth}
Nesterov, Y. (2005).
\newblock Smooth minimization of non-smooth functions.
\newblock {\em Mathematical Programming}, 103(1):127--152.

\bibitem[Nesterov, 2018]{AGD3}
Nesterov, Y. (2018).
\newblock {\em Lectures on convex optimization}, volume 137.
\newblock Springer.

\bibitem[Nguyen et~al., 2017]{SARAH}
Nguyen, L.~M., Liu, J., Scheinberg, K., and Tak{\'a}{\v{c}}, M. (2017).
\newblock {SARAH: A Novel Method for Machine Learning Problems Using Stochastic
  Recursive Gradient}.
\newblock In {\em Proceedings of the 34th International Conference on Machine
  Learning}, pages 2613--2621.

\bibitem[Nguyen et~al., 2018]{NguyenNDRST18}
Nguyen, L.~M., Nguyen, P.~H., van Dijk, M., Richt{\'{a}}rik, P., Scheinberg,
  K., and Tak{\'{a}}c, M. (2018).
\newblock {{SGD} and Hogwild! Convergence Without the Bounded Gradients
  Assumption}.
\newblock In {\em Proceedings of the 35th International Conference on Machine
  Learning}, volume~80, pages 3747--3755.

\bibitem[Pedregosa et~al., 2017]{ProxASAGA}
Pedregosa, F., Leblond, R., and Lacoste{-}Julien, S. (2017).
\newblock {Breaking the Nonsmooth Barrier: {A} Scalable Parallel Method for
  Composite Optimization}.
\newblock In {\em Advances in Neural Information Processing Systems}, pages
  56--65.

\bibitem[Recht et~al., 2011]{Hogwild}
Recht, B., R{\'{e}}, C., Wright, S.~J., and Niu, F. (2011).
\newblock {Hogwild!: A Lock-Free Approach to Parallelizing Stochastic Gradient
  Descent}.
\newblock In {\em Advances in Neural Information Processing Systems}, pages
  693--701.

\bibitem[Roux et~al., 2012]{SAG}
Roux, N.~L., Schmidt, M., and Bach, F.~R. (2012).
\newblock {A Stochastic Gradient Method with an Exponential Convergence Rate
  for Finite Training Sets}.
\newblock In {\em Advances in Neural Information Processing Systems}, pages
  2663--2671.

\bibitem[Sa et~al., 2015]{SaZORR15}
Sa, C.~D., Zhang, C., Olukotun, K., and R{\'{e}}, C. (2015).
\newblock {Taming the Wild: {A} Unified Analysis of Hogwild-Style Algorithms}.
\newblock In {\em Advances in Neural Information Processing Systems}, pages
  2674--2682.

\bibitem[Schmidt et~al., 2017]{SAG2}
Schmidt, M., Le~Roux, N., and Bach, F. (2017).
\newblock Minimizing finite sums with the stochastic average gradient.
\newblock {\em Mathematical Programming}, 162(1-2):83--112.

\bibitem[Shang et~al., 2018]{ShangLZCNY18}
Shang, F., Liu, Y., Zhou, K., Cheng, J., Ng, K. K.~W., and Yoshida, Y. (2018).
\newblock {Guaranteed Sufficient Decrease for Stochastic Variance Reduced
  Gradient Optimization}.
\newblock In {\em Proceedings of the Twenty-First International Conference on
  Artificial Intelligence and Statistics}, volume~84, pages 1027--1036.

\bibitem[Shi et~al., 2021]{shi2021ai}
Shi, Z., Loizou, N., Richt{\'a}rik, P., and Tak{\'a}{\v{c}}, M. (2021).
\newblock {AI-SARAH: Adaptive and Implicit Stochastic Recursive Gradient
  Methods}.
\newblock {\em arXiv preprint arXiv:2102.09700}.

\bibitem[Song et~al., 2020]{VRADA}
Song, C., Jiang, Y., and Ma, Y. (2020).
\newblock {Variance Reduction via Accelerated Dual Averaging for Finite-Sum
  Optimization}.
\newblock In {\em Advances in Neural Information Processing Systems},
  volume~33, pages 833--844.

\bibitem[Stich et~al., 2021]{StichMJ21}
Stich, S.~U., Mohtashami, A., and Jaggi, M. (2021).
\newblock {Critical Parameters for Scalable Distributed Learning with Large
  Batches and Asynchronous Updates}.
\newblock In {\em Proceedings of the Twenty-Fourth International Conference on
  Artificial Intelligence and Statistics}, volume 130, pages 4042--4050.

\bibitem[Tan et~al., 2016]{TanMDQ16}
Tan, C., Ma, S., Dai, Y., and Qian, Y. (2016).
\newblock {Barzilai-Borwein Step Size for Stochastic Gradient Descent}.
\newblock In {\em Advances in Neural Information Processing Systems}, pages
  685--693.

\bibitem[Taylor et~al., 2017]{taylor2017smooth}
Taylor, A.~B., Hendrickx, J.~M., and Glineur, F. (2017).
\newblock Smooth strongly convex interpolation and exact worst-case performance
  of first-order methods.
\newblock {\em Mathematical Programming}, 161(1-2):307--345.

\bibitem[Woodworth and Srebro, 2016]{woodworthS16}
Woodworth, B.~E. and Srebro, N. (2016).
\newblock {Tight Complexity Bounds for Optimizing Composite Objectives}.
\newblock In {\em Advances in Neural Information Processing Systems}, pages
  3639--3647.

\bibitem[Xiao and Zhang, 2014]{Prox_SVRG}
Xiao, L. and Zhang, T. (2014).
\newblock {A Proximal Stochastic Gradient Method with Progressive Variance
  Reduction}.
\newblock {\em SIAM Journal on Optimization}, 24(4):2057--2075.

\bibitem[Yu et~al., 2010]{yu2010feature}
Yu, H.-F., Lo, H.-Y., Hsieh, H.-P., Lou, J.-K., McKenzie, T.~G., Chou, J.-W.,
  Chung, P.-H., Ho, C.-H., Chang, C.-F., Wei, Y.-H., et~al. (2010).
\newblock {Feature engineering and classifier ensemble for KDD cup 2010}.
\newblock In {\em KDD cup}.

\bibitem[Zhou et~al., 2019]{SSNM}
Zhou, K., Ding, Q., Shang, F., Cheng, J., Li, D., and Luo, Z.-Q. (2019).
\newblock {Direct Acceleration of SAGA using Sampled Negative Momentum}.
\newblock In {\em Proceedings of the Twenty-Second International Conference on
  Artificial Intelligence and Statistics}, pages 1602--1610.

\bibitem[Zhou et~al., 2020a]{pmlr-v124-zhou20a}
Zhou, K., Jin, Y., Ding, Q., and Cheng, J. (2020a).
\newblock {Amortized Nesterov’s Momentum: A Robust Momentum and Its
  Application to Deep Learning}.
\newblock In {\em Proceedings of the 36th Conference on Uncertainty in
  Artificial Intelligence}, pages 211--220.

\bibitem[Zhou et~al., 2018]{MiG}
Zhou, K., Shang, F., and Cheng, J. (2018).
\newblock {A Simple Stochastic Variance Reduced Algorithm with Fast Convergence
  Rates}.
\newblock In {\em Proceedings of the 35th International Conference on Machine
  Learning}, pages 5980--5989.

\bibitem[Zhou et~al., 2020b]{zhou2020boosting}
Zhou, K., So, A. M.-C., and Cheng, J. (2020b).
\newblock {Boosting First-Order Methods by Shifting Objective: New Schemes with
  Faster Worst-Case Rates}.
\newblock In {\em Advances in Neural Information Processing Systems}, pages
  15405--15416.

\bibitem[Zhou et~al., 2021]{zhou2021practical}
Zhou, K., Tian, L., So, A. M.-C., and Cheng, J. (2021).
\newblock {Practical Schemes for Finding Near-Stationary Points of Convex
  Finite-Sums}.
\newblock {\em arXiv preprint arXiv:2105.12062}.

\end{thebibliography}
% Supplementary material: To improve readability, you must use a single-column format for the supplementary material.
\appendix

\section{Proof of Lemma \ref{lem:S-SVRG-VB}}
\label{app:lemma_vb}
\[
\begin{aligned}
	&\Ei{\norm{\mathcal{G}_y - \pf{y}}^2} \\
	\meq{a}{}& \Ei{\norm{\pfi{y} - \pfi{\tilde{x}} + D_{i}\pf{\tilde{x}}}^2} - \norm{\pf{y}}^2\\
	={}& \Ei{\norm{\pfi{y} - \pfi{\tilde{x}}}^2} + 2 \Ei{\innr{\pfi{y} - \pfi{\tilde{x}}, D_{i} \pf{\tilde{x}}}} + \Ei{\norm{D_{i} \pf{\tilde{x}}}^2} - \norm{\pf{y}}^2\\
	\meq{b}{}& \Ei{\norm{\pfi{y} - \pfi{\tilde{x}}}^2} + 2 \Ei{\innr{\pfi{y} - \pfi{\tilde{x}}, D \pf{\tilde{x}}}} + \Ei{\innr{D\pf{\tilde{x}}, D_{i}\pf{\tilde{x}}}} - \norm{\pf{y}}^2\\
	={}& \Ei{\norm{\pfi{y} - \pfi{\tilde{x}}}^2} + 2\innr{\pf{y}, D \pf{\tilde{x}}} - \innr{\pf{\tilde{x}},D\pf{\tilde{x}}} - \norm{\pf{y}}^2\\
	\mleq{c}{}& 2L\big(f(\tilde{x}) - f(y) - \innr{\pf{y}, \tilde{x} - y}\big) + 2\innr{\pf{y}, D \pf{\tilde{x}}} - \innr{\pf{\tilde{x}},D\pf{\tilde{x}}} - \norm{\pf{y}}^2,
\end{aligned}
\]
where $\mar{a}$ uses the unbiasedness $\Ei{\mathcal{G}_y} = \pf{y}$, $\mar{b}$ follows from that $\pfi{y} - \pfi{\tilde{x}}$ is supported on $T_{i}$ and $P_{i}^2 = P_{i}$, and $\mar{c}$ uses the interpolation condition of $f_{i}$ (see Section \ref{sec:preliminaries}).

\section{Proof of Theorem \ref{thm:SS-Acc-SVRG}}
\label{app:thm_serial}
We omit the superscript $s$ in the one-epoch analysis for clarity.
Using convexity, we have
\begin{align}
	f(y_k) - f(\xs) \leq{}& \innr{\pf{y_k}, y_k - \xs} \nonumber \\
	={}& \innr{\pf{y_k}, y_k - z_k} + \innr{\pf{y_k}, z_k - \xs}\nonumber\\
	\mleq{\star}{}& \frac{1 - \vartheta}{\vartheta}\innr{\pf{y_k}, \tilde{x}_s - y_k} - \frac{\varphi}{\vartheta}\innr{\pf{y_k},D\pf{\tilde{x}_s}} + \innr{\pf{y_k}, z_k - \xs},\label{P7}
\end{align}
where $\mar{\star}$ follows from the construction $y_{k} = \vartheta z_k + \left(1 - \vartheta\right) \tilde{x}_s - \varphi D\pf{\tilde{x}_s}$.

Denote $\mathcal{G}_{y_k} = \pfik{y_k} - \pfik{\tilde{x}_s} + D_{i_k}\pf{\tilde{x}_s}$. Based on the updating rule: $z_{k+1} = z_k - \eta \cdot \mathcal{G}_{y_k}$, it holds that
\begin{align}
	&\norm{z_{k+1} - \xs}^2 = \norm{z_k - \xs - \eta \cdot \mathcal{G}_{y_k}}^2 \nonumber\\ \Rightarrow{}&
	\innr{\mathcal{G}_{y_k}, z_k - \xs} = \frac{\eta}{2}\norm{\mathcal{G}_{y_k}}^2 + \frac{1}{2\eta} \left(\norm{z_k - \xs}^2 - \norm{z_{k+1} - \xs}^2\right) \nonumber\\
	\mra{\star}{}&\innr{\pf{y_k}, z_k - \xs} = \frac{\eta}{2}\Eik{\norm{\mathcal{G}_{y_k}}^2} + \frac{1}{2\eta} \left(\norm{z_k - \xs}^2 - \Eik{\norm{z_{k+1} - \xs}^2}\right)\nonumber\\
	\Rightarrow{}& \innr{\pf{y_k}, z_k - \xs} = \frac{\eta}{2}\Eik{\norm{\mathcal{G}_{y_k} - \pf{y_k}}^2} + \frac{\eta}{2}\norm{\pf{y_k}}^2 + \frac{1}{2\eta} \left(\norm{z_k - \xs}^2 - \Eik{\norm{z_{k+1} - \xs}^2}\right), \label{P8}
\end{align}
where $\mar{\star}$ follows from taking the expectation with respect to  sample $i_k$.

Combine \eqref{P7}, \eqref{P8} and use the variance bound in Lemma \ref{lem:S-SVRG-VB}.
\[ 
\begin{aligned}
	f(y_k) - f(\xs)
	\leq{}& \eta L\big(f(\tilde{x}_s) - f(y_k)\big) + \left(\frac{1 - \vartheta}{\vartheta}- \eta L\right)\innr{\pf{y_k}, \tilde{x}_s - y_k} + \left(\eta- \frac{\varphi}{\vartheta}\right) \innr{\pf{y_k}, D \pf{\tilde{x}_s}} \\&-\frac{\eta}{2}\innr{\pf{\tilde{x}_s},D\pf{\tilde{x}_s}} + \frac{1}{2\eta} \left(\norm{z_k - \xs}^2 - \Eik{\norm{z_{k+1} - \xs}^2}\right).
\end{aligned}
\]

Substituting the choices $\eta = \frac{1 - \vartheta}{L\vartheta}$ and $\varphi = \eta\vartheta = \frac{1 - \vartheta}{L}$ and noting that $D \succ 0$, we obtain
\[ 
\begin{aligned}
	f(y_k) - f(\xs)
	\leq{}& (1 - \vartheta)\big(f(\tilde{x}_s) - f(\xs)\big) + \frac{L\vartheta^2}{2(1 - \vartheta)} \left(\norm{z_k - \xs}^2 - \Eik{\norm{z_{k+1} - \xs}^2}\right).
\end{aligned}
\]

Summing the above inequality from $k=0$ to $m-1$ and noting that $z^{s+1}_0 = z^s_m$, we obtain
\[ 
\begin{aligned}
	\E{f(\tilde{x}_{s+1}) - f(\xs)} ={}& \frac{1}{m} \sum_{k=0}^{m-1}{\E{f(y_k) - f(\xs)}}\\
	\leq{}& (1 - \vartheta)\E{f(\tilde{x}_s) - f(\xs)} + \frac{L\vartheta^2}{2m(1-\vartheta)} \left(\E{\norm{z^s_0- \xs}^2} - \E{\norm{z^{s+1}_0 - \xs}^2}\right).
\end{aligned}
\]

Denoting $Q_s \triangleq \E{f(\tilde{x}_s) - f(\xs)}$ and $P_s \triangleq \E{\norm{z^s_0- \xs}^2}$, we can write the above as
\[ 
\begin{aligned}
	Q_{s+1}
	\leq{}& \frac{1 - \vartheta}{\vartheta}(Q_s - Q_{s+1}) + \frac{L\vartheta}{2m(1-\vartheta)} \left(P_s - P_{s+1}\right).
\end{aligned}
\]

Summing this inequality from $s = 0$ to $S-1$ and using Jensen's inequality, we have
\[
\begin{aligned}
	\E{f(x_{r+1}) - f(\xs)} \leq{}& \frac{1}{S}\sum_{s=0}^{S-1}{Q_{s+1}} \\
	\leq{}& \frac{1 - \vartheta}{\vartheta S}(Q_0 - Q_S) + \frac{L\vartheta}{2m(1-\vartheta)S} \left(P_0 - P_S\right) \\
	\mleq{\star}{}& \frac{1 - \vartheta}{\vartheta S} \E{f(x_r) - f(\xs)} + \frac{L\vartheta}{2m(1-\vartheta)S} \E{\norm{x_r - \xs}^2},
\end{aligned}
\]
where $\mar{\star}$ follows from $\tilde{x}_0 = z^0_0 = x_r$.

Using the $\mu$-strong convexity at $(x_r, \xs)$ (or quadratic growth), i.e., $f(x_r) - f(\xs) \geq \frac{\mu}{2}\norm{x_r - \xs}^2$, we arrive at
\[
\E{f(x_{r+1}) - f(\xs)} \leq \left(\frac{1 - \vartheta}{\vartheta S} + \frac{\kappa\vartheta}{m(1-\vartheta)S} \right) \E{f(x_r) - f(\xs)}. 
\]

Choosing $S = \left\lceil \omega \cdot \left(\frac{1 - \vartheta}{\vartheta} + \frac{\kappa\vartheta}{m(1-\vartheta)}\right) \right\rceil$, we have $\E{f(x_{r+1}) - f(\xs)} \leq \frac{1}{\omega} \cdot \E{f(x_r) - f(\xs)}$. Thus, for any accuracy $\epsilon > 0$ and any constant $\omega > 1$, to guarantee that the output satisfies $\E{f(x_R)} - f(\xs) \leq \epsilon$, we need to perform totally $R = O\left(\log{\frac{f(x_0) - f(\xs)}{\epsilon}}\right)$ restarts. Note that the total oracle complexity of Algorithm \ref{alg:SS-Acc-SVRG} is 
$\text{\#grad} = R\cdot S \cdot (n+2m)$. Setting $m = \Theta(n)$ and minimizing $S$ with respect to $\vartheta$, we obtain the optimal choice $\vartheta = \frac{\sqrt{m}}{\sqrt{\kappa} + \sqrt{m}}$. In this case, we have $S =  \left\lceil 2\omega\sqrt{\frac{\kappa}{m}}\right\rceil = O\left(\max{\left\{1, \sqrt{\frac{\kappa}{n}}\right\}}\right)$. Finally, the total oracle complexity is
\[
\text{\#grad} = O\left(\max{\left\{n, \sqrt{\kappa n}\right\}}\log{\frac{f(x_0) - f(\xs)}{\epsilon}}\right).
\]

\newpage
\section{Proof of Theorem \ref{thm:AS-Acc-SVRG}}
\label{app:thm_async}

We use the following lemma from \citep{ASAGA_JMLR} to bound the term in asynchrony. Although the gradient is evaluated at a different point, the same proof works in our case, and thus is omitted here.
\begin{lemmas}[Inequality (63) in \citep{ASAGA_JMLR}]
	\label{lem:async_term_bound}
	We have the following bound for the iterates $z_k$ (defined at \eqref{virtual_iterate}) and $\hat{z}_k$ (defined at \eqref{def_inconsistency}) in one epoch of Algorithm \ref{alg:AS-Acc-SVRG}:
	\[
	\E{\innr{\mathcal{G}_{\hat{y}_k}, \hat{z}_k - z_k}}\leq \frac{\sqrt{\Delta}\eta}{2}\left(\sum_{j = (k-\tau)_+}^{k-1} {\E{\norm{\mathcal{G}_{\hat{y}_j}}^2} + \tau\E{\norm{\mathcal{G}_{\hat{y}_k}}^2}}\right).
	\]
\end{lemmas}

We start with analyzing one epoch of virtual iterates defined at \eqref{virtual_iterate}: $z_{k+1} = z_k - \eta\cdot\mathcal{G}_{\hat{y}_k}$. The expectation is first taken conditioned on the previous epochs.
\begin{align}
	&\norm{z_{k+1} - \xs}^2 = \norm{z_k - \xs - \eta \cdot \mathcal{G}_{\hat{y}_k}}^2 \nonumber\\ \Rightarrow{}&
	\innr{\mathcal{G}_{\hat{y}_k}, \hat{z}_k - \xs} + \innr{\mathcal{G}_{\hat{y}_k}, z_k - \hat{z}_k} = \frac{\eta}{2}\norm{\mathcal{G}_{\hat{y}_k}}^2 + \frac{1}{2\eta} \left(\norm{z_k - \xs}^2 - \norm{z_{k+1} - \xs}^2\right) \nonumber\\
	\mra{\star}{}&\E{\innr{\pf{\hat{y}_k}, \hat{z}_k - \xs}} = \frac{\eta}{2}\E{\norm{\mathcal{G}_{\hat{y}_k}}^2} + \E{\innr{\mathcal{G}_{\hat{y}_k}, \hat{z}_k - z_k}} + \frac{1}{2\eta} \left(\E{\norm{z_k- \xs}^2} - \E{\norm{z_{k+1 } - \xs}^2}\right), \label{P9}
\end{align}
where $\mar{\star}$ follows from taking the expectation and the unbiasedness assumption $\E{\mathcal{G}_{\hat{y}_k} | \hat{z}_k} = \pf{\hat{y}_k}$. 

Using the convexity at $\xs$ and the inconsistent $\hat{y}_k$ (ordered at \eqref{def_inconsistency_y}), we have 
\[
\begin{aligned}
	f(\hat{y}_k) - f(\xs) \leq{}& \innr{\pf{\hat{y}_k}, \hat{y}_k - \xs}  \\
	={}& \innr{\pf{\hat{y}_k}, \hat{y}_k - \hat{z}_k} + \innr{\pf{\hat{y}_k}, \hat{z}_k - \xs}
	\\
	\mleq{\star}{}& \frac{1 - \vartheta}{\vartheta}\innr{\pf{\hat{y}_k}, \tilde{x}_s - \hat{y}_k} - \frac{\varphi}{\vartheta}\innr{\pf{\hat{y}_k},D\pf{\tilde{x}_s}} + \innr{\pf{\hat{y}_k}, \hat{z}_k - \xs}.
\end{aligned}
\]
where $\mar{\star}$ follows from the construction $\hat{y}_{k} = \vartheta \hat{z}_k + \left(1 - \vartheta\right) \tilde{x}_s - \varphi D\pf{\tilde{x}_s}$.

After taking the expectation and combining with \eqref{P9}, we obtain
\[
\begin{aligned}
	\E{f(\hat{y}_k)} - f(\xs) \leq{}& \frac{1 - \vartheta}{\vartheta}\E{\innr{\pf{\hat{y}_k}, \tilde{x}_s - \hat{y}_k}} - \frac{\varphi}{\vartheta}\E{\innr{\pf{\hat{y}_k},D\pf{\tilde{x}_s}}}+\frac{\eta}{2}\E{\norm{\mathcal{G}_{\hat{y}_k}}^2} \\& + \E{\innr{\mathcal{G}_{\hat{y}_k}, \hat{z}_k - z_k}} + \frac{1}{2\eta} \left(\E{\norm{z_k - \xs}^2} - \E{\norm{z_{k+1} - \xs}^2}\right).
\end{aligned}
\]

Summing this inequality from $k = 0$ to $m-1$, we have
\[
\begin{aligned}
	\sum_{k=0}^{m-1}{\E{f(\hat{y}_k) - f(\xs)}} \leq{}& \sum_{k=0}^{m-1}{\left(\frac{1 - \vartheta}{\vartheta}\E{\innr{\pf{\hat{y}_k}, \tilde{x}_s - \hat{y}_k}} - \frac{\varphi}{\vartheta}\E{\innr{\pf{\hat{y}_k},D\pf{\tilde{x}_s}}} +\frac{\eta}{2}\E{\norm{\mathcal{G}_{\hat{y}_k}}^2}\right)} \\&+ \sum_{k=0}^{m-1}{\E{\innr{\mathcal{G}_{\hat{y}_k}, \hat{z}_k - z_k}}} + \frac{1}{2\eta} \left(\norm{z_0 - \xs}^2 - \E{\norm{z_m - \xs}^2}\right).
\end{aligned}
\]

Invoke Lemma \ref{lem:async_term_bound} to bound the asynchronous perturbation.
\[
\begin{aligned}
	&\sum_{k=0}^{m-1}{\E{f(\hat{y}_k) - f(\xs)}} \\ \leq{}& \sum_{k=0}^{m-1}{\left(\frac{1 - \vartheta}{\vartheta}\E{\innr{\pf{\hat{y}_k}, \tilde{x}_s - \hat{y}_k}} - \frac{\varphi}{\vartheta}\E{\innr{\pf{\hat{y}_k},D\pf{\tilde{x}_s}}} +\frac{\eta}{2}\E{\norm{\mathcal{G}_{\hat{y}_k}}^2}\right)} \\&+ \frac{\sqrt{\Delta}\eta}{2}\sum_{k=0}^{m-1}{\left(\sum_{j = (k-\tau)_+}^{k-1} {\E{\norm{\mathcal{G}_{\hat{y}_j}}^2} + \tau\E{\norm{\mathcal{G}_{\hat{y}_k}}^2}}\right)} + \frac{1}{2\eta} \left(\norm{z_0 - \xs}^2 - \E{\norm{z_m - \xs}^2}\right)\\
	={}& \sum_{k=0}^{m-1}{\left(\frac{1 - \vartheta}{\vartheta}\E{\innr{\pf{\hat{y}_k}, \tilde{x}_s - \hat{y}_k}} - \frac{\varphi}{\vartheta}\E{\innr{\pf{\hat{y}_k},D\pf{\tilde{x}_s}}} +\frac{\eta(1 + \sqrt{\Delta}\tau)}{2}\E{\norm{\mathcal{G}_{\hat{y}_k}}^2}\right)} \\&+ \frac{\sqrt{\Delta}\eta}{2}\sum_{k=0}^{m-1}{\sum_{j = (k-\tau)_+}^{k-1} {\E{\norm{\mathcal{G}_{\hat{y}_j}}^2}}} + \frac{1}{2\eta} \left(\norm{z_0 - \xs}^2 - \E{\norm{z_m - \xs}^2}\right) \\
	\leq{}& \sum_{k=0}^{m-1}{\left(\frac{1 - \vartheta}{\vartheta}\E{\innr{\pf{\hat{y}_k}, \tilde{x}_s - \hat{y}_k}} - \frac{\varphi}{\vartheta}\E{\innr{\pf{\hat{y}_k},D\pf{\tilde{x}_s}}} +\frac{\eta(1 + \sqrt{\Delta}\tau)}{2}\E{\norm{\mathcal{G}_{\hat{y}_k}}^2}\right)} \\&+ \frac{\sqrt{\Delta}\eta\tau}{2}\sum_{k=0}^{m-1}{\E{\norm{\mathcal{G}_{\hat{y}_k}}^2}} + \frac{1}{2\eta} \left(\norm{z_0 - \xs}^2 - \E{\norm{z_m - \xs}^2}\right)\\
	={}& \sum_{k=0}^{m-1}{\left(\frac{1 - \vartheta}{\vartheta}\E{\innr{\pf{\hat{y}_k}, \tilde{x}_s - \hat{y}_k}} - \frac{\varphi}{\vartheta}\E{\innr{\pf{\hat{y}_k},D\pf{\tilde{x}_s}}} +\frac{\eta(1 + 2\sqrt{\Delta}\tau)}{2}\E{\norm{\mathcal{G}_{\hat{y}_k}}^2}\right)} \\&+ \frac{1}{2\eta} \left(\norm{z_0 - \xs}^2 - \E{\norm{z_m - \xs}^2}\right). 
\end{aligned}
\]

For some $\widetilde{\tau} \geq \tau$, we choose $\eta = \frac{(1 - \vartheta)}{L\vartheta(1 + 2\sqrt{\Delta}\widetilde{\tau})}$, and then
\[
\begin{aligned}
	\sum_{k=0}^{m-1}{\E{f(\hat{y}_k) - f(\xs)}} \leq{}&\sum_{k=0}^{m-1}{\left(\frac{1 - \vartheta}{\vartheta}\E{\innr{\pf{\hat{y}_k}, \tilde{x}_s - \hat{y}_k}} - \frac{\varphi}{\vartheta}\E{\innr{\pf{\hat{y}_k},D\pf{\tilde{x}_s}}} +\frac{1 - \vartheta}{2L\vartheta}\E{\norm{\mathcal{G}_{\hat{y}_k}}^2}\right)} \\&+ \frac{L\vartheta(1+2\sqrt{\Delta}\widetilde{\tau})}{2(1-\vartheta)} \left(\norm{z_0 - \xs}^2 - \E{\norm{z_m - \xs}^2}\right).
\end{aligned}
\]
Note that $\E{\norm{\mathcal{G}_{\hat{y}_k}}^2} = \E{\norm{\mathcal{G}_{\hat{y}_k} - \pf{\hat{y}_k}}^2} + \E{\norm{\pf{\hat{y}_k}}^2}$ due to the unbiasedness assumption. Using Lemma \ref{lem:S-SVRG-VB}, we can conclude that
\[
\begin{aligned}
	\sum_{k=0}^{m-1}{\E{f(\hat{y}_k) - f(\xs)}}\leq{}&\sum_{k=0}^{m-1}{\left(\left(\frac{1 - \vartheta}{L\vartheta} - \frac{\varphi}{\vartheta}\right)\E{\innr{\pf{\hat{y}_k},D\pf{\tilde{x}_s}}} +\frac{1 - \vartheta}{\vartheta}\big(f(\tilde{x}_s) - \E{f(\hat{y}_k)}\big)\right)} \\&+ \frac{L\vartheta(1+2\sqrt{\Delta}\widetilde{\tau})}{2(1-\vartheta)} \left(\norm{z_0 - \xs}^2 - \E{\norm{z_m - \xs}^2}\right).
\end{aligned}
\]

Choosing $\varphi = \frac{1 - \vartheta}{L}$ and dividing both sides by $m$, we can arrange this inequality as
\[
\begin{aligned}
	\frac{1}{m}\sum_{k=0}^{m-1}{\E{f(\hat{y}_k) - f(\xs)}}
	\leq (1-\vartheta)\big(f(\tilde{x}_s) - f(\xs)\big) + \frac{L\vartheta^2(1+2\sqrt{\Delta}\widetilde{\tau})}{2m(1 - \vartheta)} \left(\norm{z_0 - \xs}^2 - \E{\norm{z_m - \xs}^2}\right).
\end{aligned}
\]

Since $\tilde{x}_{s+1}$ is chosen uniformly at random from $\{\hat{y}_0, \ldots, \hat{y}_{m-1}\}$ and that $z^{s+1}_0 = z^s_m$ (the first and the last virtual iterates exist in the shared memory, and $z$ is unchanged after each epoch in Algorithm \ref{alg:AS-Acc-SVRG}), the following holds
\[
\begin{gathered}
	\E{f(\tilde{x}_{s+1}) - f(\xs)}
	\leq (1-\vartheta)\big(f(\tilde{x}_s) - f(\xs)\big) + \frac{L\vartheta^2(1+2\sqrt{\Delta}\widetilde{\tau})}{2m(1-\vartheta)} \left(\norm{z^s_0 - \xs}^2 - \E{\norm{z^{s+1}_0 - \xs}^2}\right).
\end{gathered}
\]

For the sake of clarity, we denote $Q_s \triangleq \E{f(\tilde{x}_s) - f(\xs)}$ and $P_s \triangleq \E{\norm{z^s_0 - \xs}^2}$. Then, it holds that
\[
Q_{s+1}
\leq \frac{1-\vartheta}{\vartheta} (Q_s - Q_{s+1}) + \frac{L\vartheta(1+2\sqrt{\Delta}\widetilde{\tau})}{2m (1 - \vartheta)} \left(P_s - P_{s+1}\right).
\]

Summing this inequality from $s = 0$ to $S-1$ and using Jensen's inequality, we have
\[
\begin{aligned}
	\E{f(x_{r+1})} - f(\xs) \leq{}& \frac{1}{S} \sum_{s=0}^{S-1} {Q_{s+1}} \\
	\leq{}& \frac{1-\vartheta}{\vartheta S} (Q_0 - Q_S) + \frac{L\vartheta(1+2\sqrt{\Delta}\widetilde{\tau})}{2m(1-\vartheta)S} \left(P_0 - P_S\right) \\
	\leq{}& \frac{1-\vartheta}{\vartheta S} Q_0 + \frac{L\vartheta(1+2\sqrt{\Delta}\widetilde{\tau})}{2m(1-\vartheta)S} P_0.
\end{aligned}
\]

Using $\mu$-strong convexity at $(x_r, \xs)$ (or quadratic growth), we have $\norm{x_r - \xs}^2 \leq \frac{2}{\mu} \big(f(x_r) - f(\xs)\big) \Leftrightarrow P_0 \leq \frac{2}{\mu} Q_0$ and thus
\[
\begin{aligned}
	\E{f(x_{r+1})} - f(\xs)
	\leq \left(\frac{1-\vartheta}{\vartheta S} + \frac{\kappa\vartheta(1+2\sqrt{\Delta}\widetilde{\tau})}{m(1-\vartheta)S} \right)\E{f(x_r) - f(\xs)},
\end{aligned}
\]

Letting $S = \left\lceil \omega\cdot\left( \frac{1-\vartheta}{\vartheta} + \frac{\kappa\vartheta(1+2\sqrt{\Delta}\widetilde{\tau})}{m(1-\vartheta)}\right) \right\rceil$, we have $\E{f(x_{r+1})} - f(\xs) \leq \frac{1}{\omega} \cdot \E{f(x_r) - f(\xs)}$. Then, since $\omega > 1$ is a constant, to achieve an $\epsilon$-additive error, we need to restart totally $R=O\left(\log{\frac{f(x_0) - f(\xs)}{\epsilon}}\right)$ times. Note that the total oracle complexity of Algorithm \ref{alg:AS-Acc-SVRG} is $\text{\#grad} = R\cdot S\cdot (n+2m)$. Setting $m=\Theta(n)$ and choosing $\vartheta$ that minimizes $S$, we obtain the optimal choice $\vartheta = \frac{\sqrt{m}}{\sqrt{\kappa(1+2\sqrt{\Delta}\widetilde{\tau})} + \sqrt{m}}$, which leads to $S = \left\lceil 2\omega\sqrt{\frac{\kappa}{m}(1+2\sqrt{\Delta}\widetilde{\tau})} \right\rceil = O\left(\max{\left\{1, \sqrt{\frac{\kappa}{n}(1+2\sqrt{\Delta}\widetilde{\tau})} \right\}}\right)$. Finally, the total oracle complexity is 
\[
\text{\#grad} = O\left( \max{\left\{n, \sqrt{\kappa n(1+2\sqrt{\Delta}\widetilde{\tau})} \right\}}\log{\frac{f(x_0) - f(\xs)}{\epsilon}}\right).
\]
\section{The Effect of the Constant $\omega$ (the Restart Frequency)}
\label{app:restart}

\begin{figure}[H]
	\begin{center}
		\includegraphics[width=0.65\textwidth]{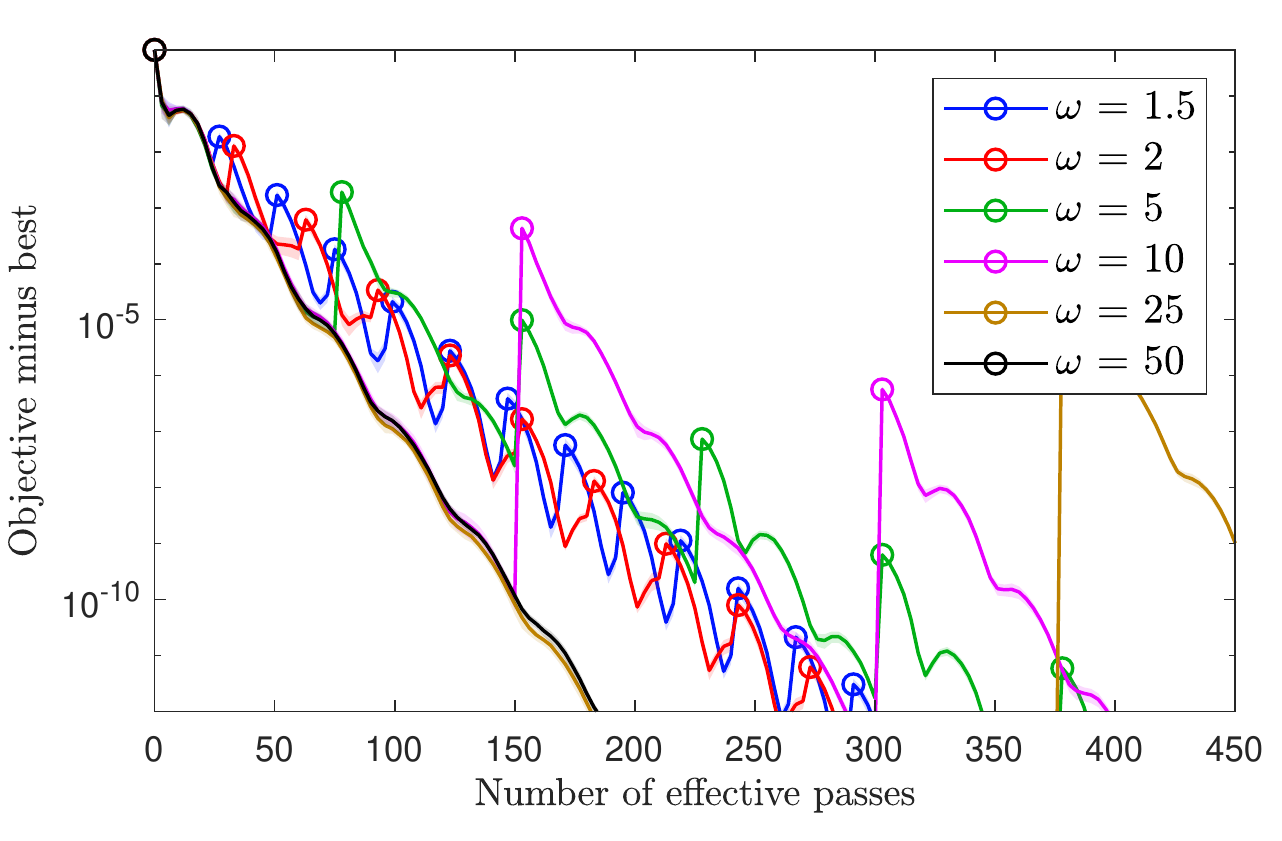}
		\caption{The practical effect of $\omega$. \textsf{RCV1.train}, run $10$ seeds. The circle marks the restarting points, i.e., $\{x_r\}$. Shaded bands indicate $\pm 1$ standard deviation.}
		\label{fig:effect_of_omega}
	\end{center}
\end{figure}
We numerically evaluate the effect of the constant $\omega$ in Figure \ref{fig:effect_of_omega}. We choose $\omega$ from $1.5$ (frequent restart) to $50$ (not restart in this task). The results on \textsf{News20} and \textsf{KDD2010.small} datasets are basically identical, and thus are omitted here. As we can see, unfortunately, the restarts only deteriorate the performance. An intuitive explanation is that the restart strategy is more conservative as it periodically retracts the point. In theory, the explicit dependence on $\omega$ (in the serial case) is 
\[
\frac{1}{\log{\omega}}\left\lceil 2\omega\sqrt{\frac{\kappa}{m}}\right\rceil,
\]
which suggests that if $\omega$ is large, the convergence on the restarting points $\{x_r\}$ will be slower. This is observed in Figure \ref{fig:effect_of_omega}. However, what our current theory cannot explain is the superior performance of not performing restart, and we are not aware a situation where the ``aggressiveness'' of no restart would hurt the convergence. Optimally tuning $\omega$ in the complexity will only leads to a small $\omega$ that closes to $1$. Further investigation is needed for the convergence of Algorithms \ref{alg:SS-Acc-SVRG} and \ref{alg:AS-Acc-SVRG} without restart (we have strong numerical evidence in Appendix \ref{app:justify} that without restart, they are still optimal).

\newpage
\section{Justifying the $\sqrt{\kappa}$ Dependence}
\label{app:justify}

\begin{figure}[H]
	\begin{center}
		\begin{subfigure}{0.325\linewidth}
			\includegraphics[width=\linewidth]{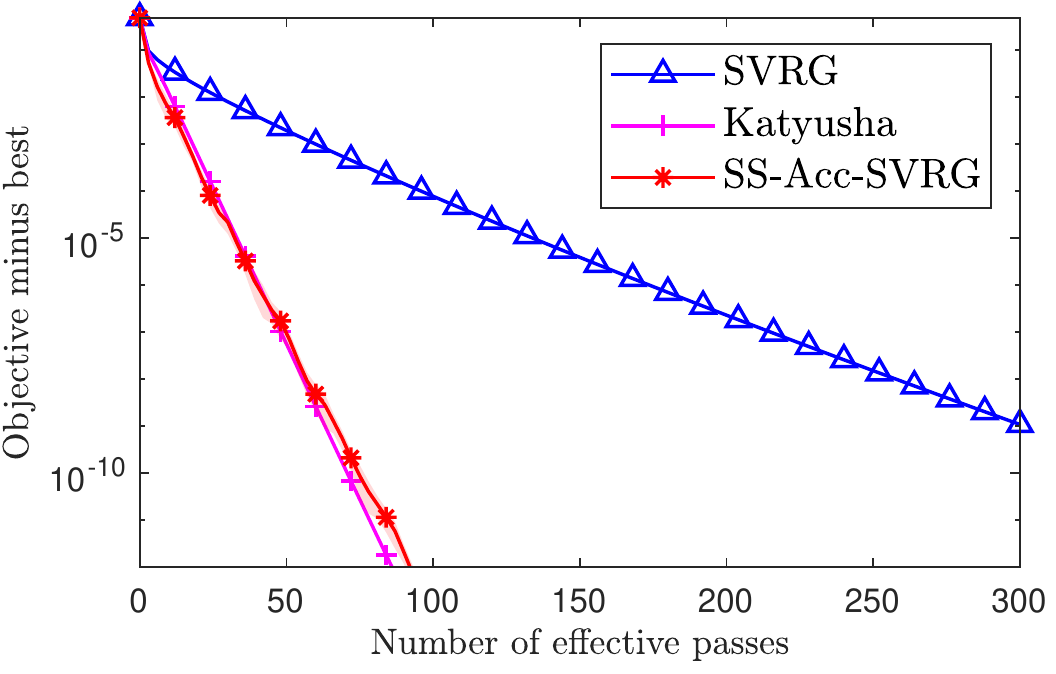}
			\caption{\textsf{KDD2010.small}, $\mu=10^{-6}$.}
			\label{fig:Grad_kdd10_small_6}
		\end{subfigure}
		\begin{subfigure}{0.325\linewidth}
			\includegraphics[width=\linewidth]{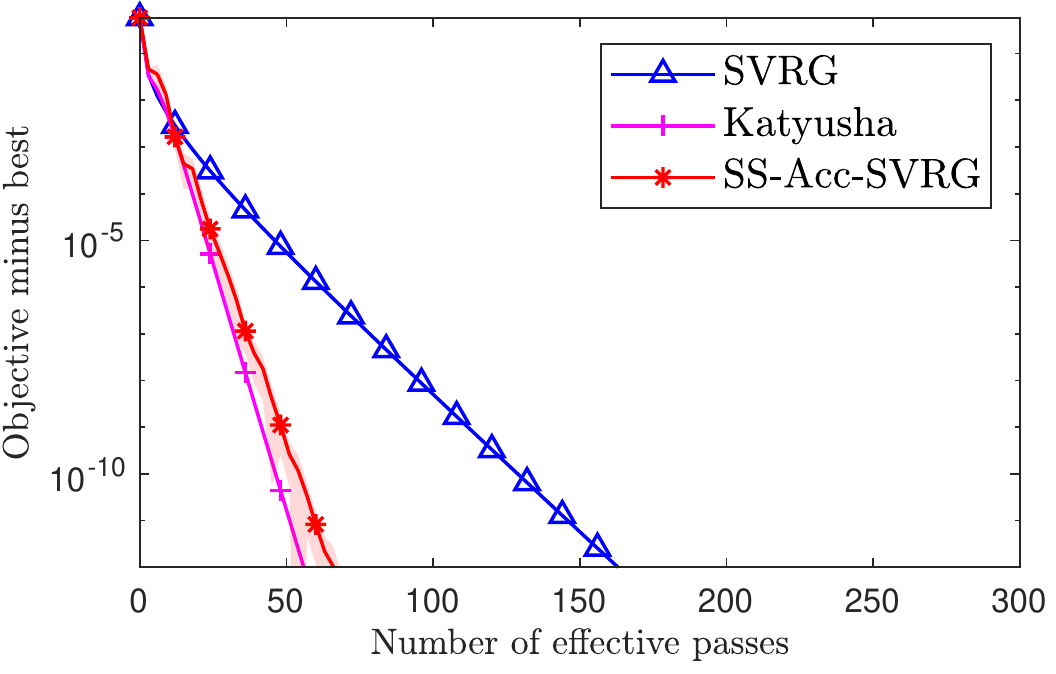}
			\caption{\textsf{RCV1.train}, $\mu=10^{-5}$.}
			\label{fig:Grad_rcv1_train_5}
		\end{subfigure} 
		\begin{subfigure}{0.325\linewidth}
			\includegraphics[width=\linewidth]{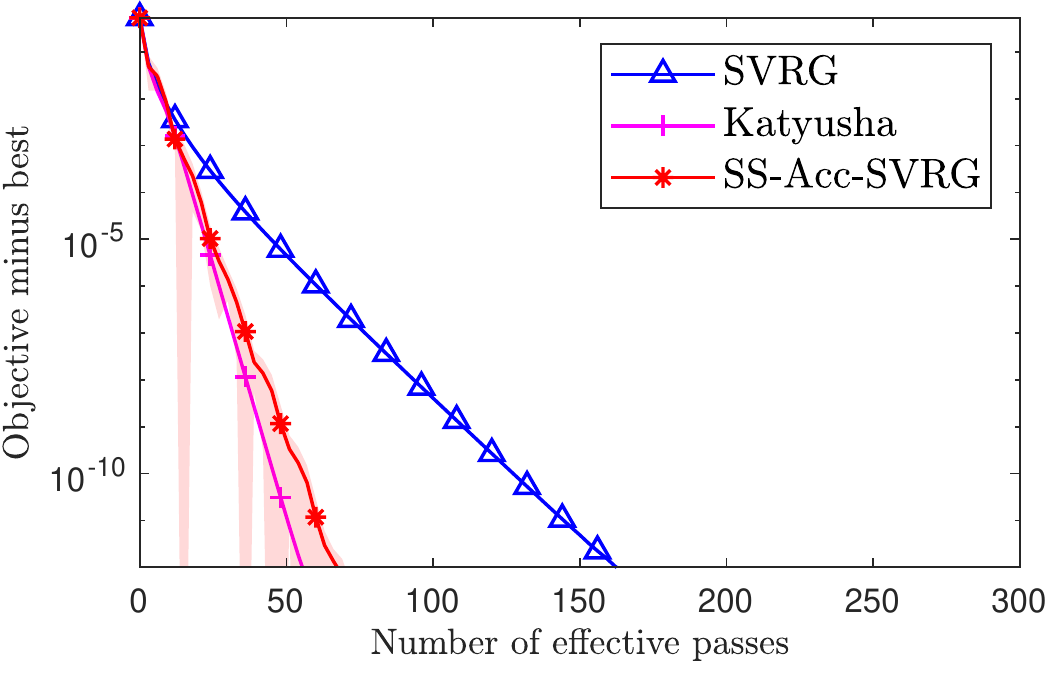}
			\caption{\textsf{News20}, $\mu=10^{-5}$.}
			\label{fig:Grad_news20_5}
		\end{subfigure} \\
		\begin{subfigure}{0.325\linewidth}
			\includegraphics[width=\linewidth]{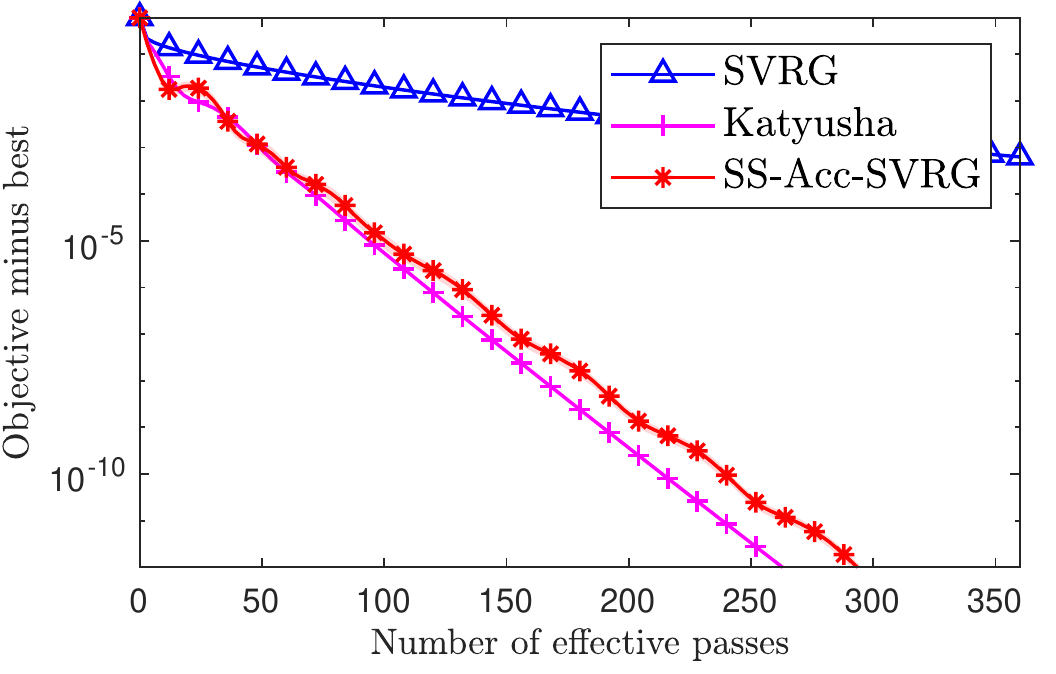}
			\caption{\textsf{KDD2010.small}, $\mu=10^{-7}$.}
			\label{fig:Grad_kdd10_small_7}
		\end{subfigure}
		\begin{subfigure}{0.325\linewidth}
			\includegraphics[width=\linewidth]{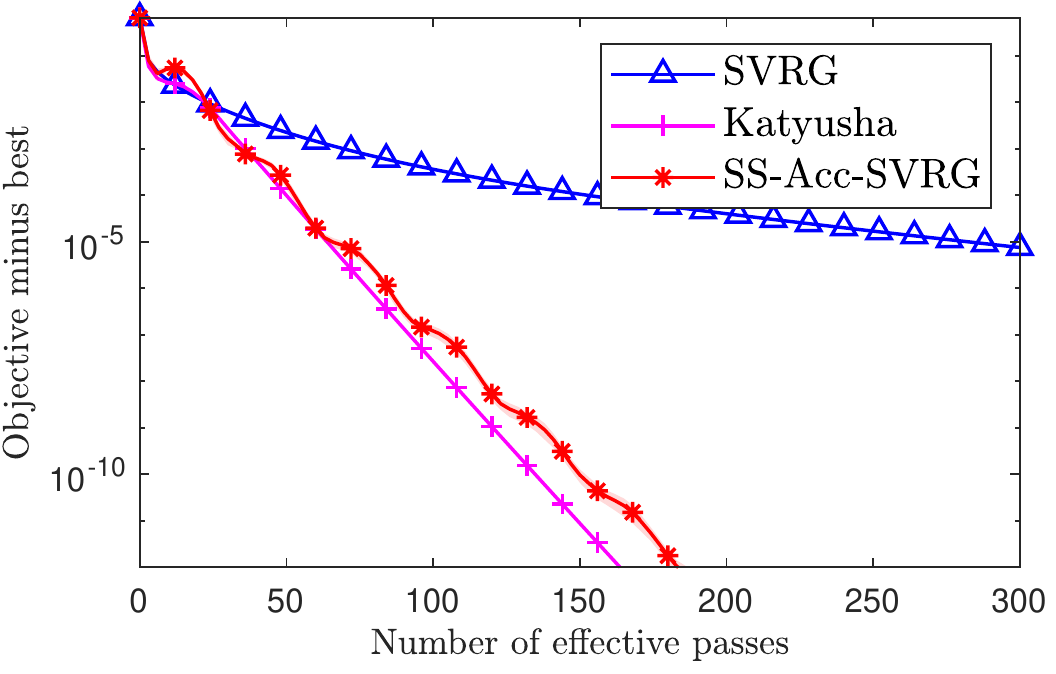}
			\caption{\textsf{RCV1.train}, $\mu=10^{-6}$.}
			\label{fig:Grad_rcv1_train_6}
		\end{subfigure} 
		\begin{subfigure}{0.325\linewidth}
			\includegraphics[width=\linewidth]{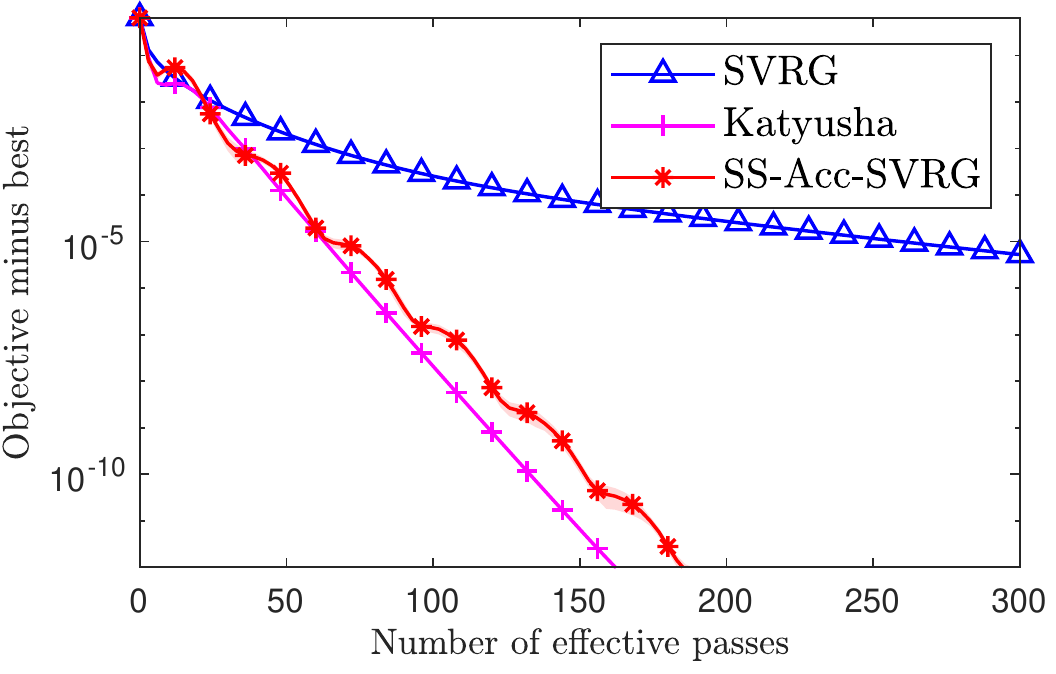}
			\caption{\textsf{News20}, $\mu=10^{-6}$.}
			\label{fig:Grad_news20_6}
		\end{subfigure}
		\caption{Justifying the $\sqrt{\kappa}$ dependence. Run $10$ seeds. Shaded bands indicate $\pm 1$ standard deviation.}
		\label{fig:Grad}
	\end{center}
\end{figure}

We choose $\mu$ to be $10$-times larger (Figure \ref{fig:Grad_kdd10_small_6} to \ref{fig:Grad_news20_5}) than the ones specified in Table \ref{data-table} (Figure \ref{fig:Grad_kdd10_small_7} to \ref{fig:Grad_news20_6}) to verify the $\sqrt{\kappa}$ dependence. In this case, $\kappa$ is $10$-times smaller and we expect the accelerated methods to be around $3$-times faster in terms of the number of data passes. This has been observed across all the datasets for Katyusha and SS-Acc-SVRG, which verifies the accelerated rate.

\section{Sanity Check for Our Implementation} 
\label{app:sanity}

\begin{figure}[H]
	\begin{center}
		\includegraphics[width=0.5\textwidth]{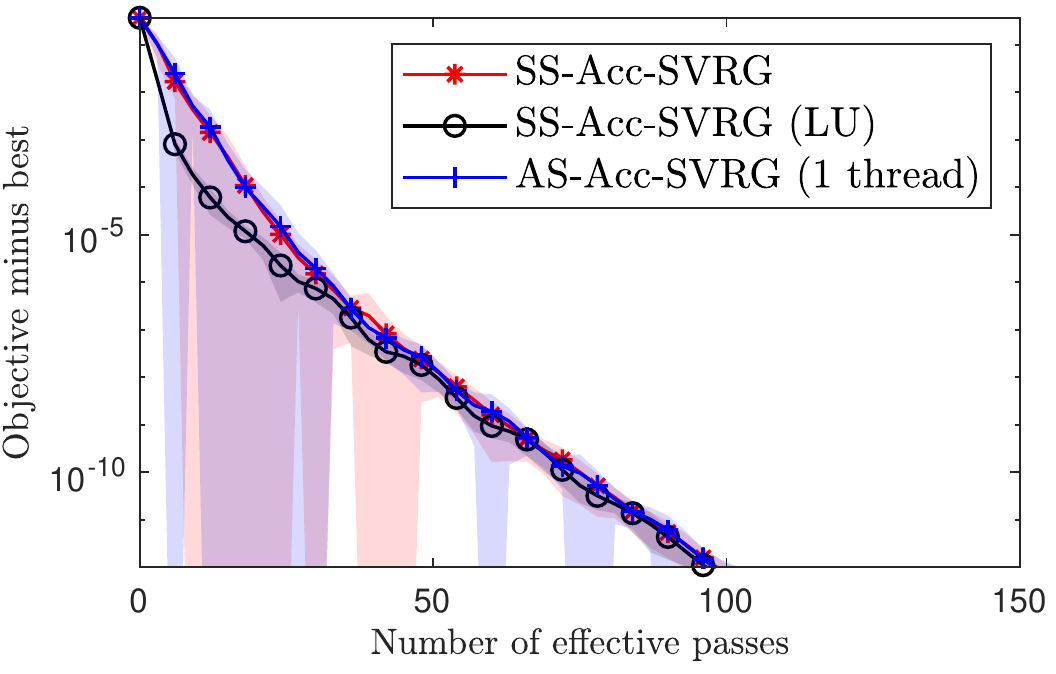}
		\caption{Sanity check. \textsf{a9a.dense}, $\mu=10^{-6}$. Run $20$ seeds. Shaded bands indicate $\pm 1$ standard deviation.}
		\label{fig:sanity_check}
	\end{center}
\end{figure}

If the dataset is fully dense, then the following three schemes should be equivalent: SS-Acc-SVRG, SS-Acc-SVRG with lagged update implementation and AS-Acc-SVRG with a single thread. We construct a fully dense dataset by adding a small positive number to each elements in \textsf{a9a} dataset \citep{LIBSVM} ($n = \SI{32561}{}, d = 123$), and we call it \textsf{a9a.dense} dataset. The sanity check is then provided in Figure \ref{fig:sanity_check}. SS-Acc-SVRG (LU) shows slightly different convergence because we used an averaged snapshot instead of a random one in its implementation. Implementing lagged update technique with a random snapshot is quite tricky, since the last seen iteration might be in the previous epochs in this case. We tested SS-Acc-SVRG with an averaged snapshot and it shows almost identical convergence as SS-Acc-SVRG (LU).

\section{Experimental Setup}
\label{app:exp}

All the methods are implemented in \CC.

\paragraph{Serial setup} We ran serial experiments on an HP Z440 machine with a single Intel Xeon E5-1630v4 with 3.70GHz cores, 16GB RAM, Ubuntu 18.04 LTS with GCC 9.4.0, MATLAB R2021a. We fixed the epoch length $m=2n$ and chose the default parameters: SS-Acc-SVRG follows Theorem \ref{alg:SS-Acc-SVRG}; Katyusha uses $\tau_2 = \frac{1}{2}, \tau_1 = \sqrt{\frac{m}{3\kappa}}, \alpha = \frac{1}{3\tau_1 L}$ \citep{allen2017katyusha}; SVRG uses $\eta = \frac{1}{4L}$.

\paragraph{Asynchronous setup} We ran asynchronous experiments on a Dell PowerEdge R840 HPC with four Intel Xeon Platinum 8160 CPU @ 2.10GHz each with 24 cores, 768GB RAM, CentOS 7.9.2009 with GCC 9.3.1, MATLAB R2021a. We implemented the original version of KroMagnon in \citep{ManiaPPRRJ17}. We fixed the epoch length $m=2n$ for KroMagnon and AS-Acc-SVRG. AS-Acc-SVRG used the same parameters as in the serial case for all datasets. The step sizes of KroMagnon and ASAGA were chosen as follows: On \textsf{RCV1.full}, KroMagnon uses $\frac{1}{L}$, ASAGA uses $\frac{1}{1.3L}$; On \textsf{Avazu-site.train}, KroMagnon uses $\frac{1}{2L}$, ASAGA uses $\frac{1}{2L}$; On \textsf{KDD2010}, KroMagnon uses $\frac{1}{2L}$, ASAGA uses $\frac{1}{3L}$. Due to the scale of the tasks, we cannot do fine-grained grid search for the step sizes of KroMagnon and ASAGA. They were chosen to ensure a stable and consistent performance in the 20-thread experiments.

\end{document}